\documentclass[11pt,a4paper]{article}
\title{Goal-Oriented A Posteriori Error Estimation for the Biharmonic Problem Based on Equilibrated Moment Tensor}

\author{Gouranga Mallik\footnote{Department of Mathematics, Indian Institute of Science, Bangalore 560012, India. Email. gourangam@iisc.ac.in}}

 \usepackage{amsmath,amsthm,amssymb,enumerate}
\usepackage{gensymb}
\usepackage[square,sort&compress,comma,numbers]{natbib}

\usepackage[sc]{mathpazo}
\usepackage{multirow} 

\usepackage{multicol}

\usepackage{newtxtext,newtxmath}
\usepackage{subfig}
\usepackage{graphicx}
\usepackage{epstopdf}
\usepackage{hyperref}
\usepackage{cancel}
\usepackage[top=30mm,bottom=30mm,foot=3mm,footskip=10mm,left=25mm,right=20mm]{geometry}
\usepackage{tikz}
\usepackage{caption}
\usetikzlibrary{shapes,calc}
\usepackage{verbatim}
\usepackage{mathrsfs}
\usepackage{algorithm}
\usepackage{accents}

\newtheorem{theorem}{Theorem}[section]
\newtheorem{lemma}[theorem]{Lemma}
\newtheorem{proof of lemma}[theorem]{Proof of Lemma}

\newtheorem{corollary}[theorem]{Corollary}

\theoremstyle{definition}
\newtheorem{definition}[theorem]{Definition}

\newtheorem{remark}[theorem]{Remark}
\numberwithin{equation}{section}

\newcommand{\bfn}{\boldsymbol{n}_e}
\newcommand{\bn}{\boldsymbol{n}}

\newcommand{\dx}{{\rm\,dx}}

\newcommand{\dss}{{\rm\,ds}}

\newcommand{\alg}[1]{\mathsf{#1}}

\newcommand{\tu}{\tilde{u}}
\newcommand{\tuh}{\tilde{u}_h}
\newcommand{\tf}{\tilde{f}}

\newcommand{\ts}{\tilde{s}}

\newcommand{\tsh}{\tilde{s}_h}

\newcommand{\sik}{\sum_{K\in\cT_h}\int_K}

\newcommand{\tbfsigh}{\tilde{\boldsymbol \sigma}_h}

\newcommand{\eqflux}{\underline{\underline{\boldsymbol \sigma}}_h^{\rm eq}}
\newcommand{\teqflux}{\underline{\underline{\tilde{\boldsymbol \sigma}}}_h^{\rm eq}}

\newcommand{\symmat}{[L^2(\Omega)]^{2\times 2}_{\rm sym}}
\newcommand{\hto}{H^2_0(\Omega)}

\newcommand{\eqfnn}{\sigma_{h,nn}^{\rm eq}}
\newcommand{\eqfnt}{\sigma_{h,n\tau}^{\rm eq}}
\newcommand{\teqfnn}{\tilde{\sigma}_{h,nn}^{\rm eq}}
\newcommand{\teqfnt}{\tilde{\sigma}_{h,n\tau}^{\rm eq}}

\newcommand{\mA}{\mathbb A}

\newcommand{\mR}{\mathbb R}

\newcommand{\mP}{\mathbb P}

\newcommand{\cE}{\mathcal E}

\newcommand{\cT}{\mathcal T}

\newcommand{\cN}{\mathcal N}

\newcommand{\cM}{\mathcal M}

\newcommand{\cR}{\mathcal R}

\newcommand{\Ctr}{C_{\text{tr}}}
\newcommand{\Ctrc}{C_{\text{tr},\text{c}}}

\newcommand{\e}{e}

\newcommand{\half}{\frac{1}{2}}
\newcommand{\Div}{{\rm div\,}}
\newcommand{\Hdiv}{{\bf H}({\rm div},\Omega)}
\newcommand{\Hdivsq}{{\bf H}({\rm div}^2,\Omega)}
\newcommand{\lt}{L^2(\Omega)}
\newcommand{\curl}{{\bf curl}}
\newcommand{\fl}{\quad \forall\:}
\newcommand{\ip}{{\rm IP}}
\newcommand{\dg}{{\rm dG}}

\newcommand{\Mh}{\underline{\underline{\bf M}}_h^{\rm eq}}
\newcommand{\tauh}{\underline{\underline{\boldsymbol{\tau}}}_h}
\newcommand{\psibarhat}{\underline{\widehat{\boldsymbol{\psi}}}}
\newcommand{\psibarheq}{\underline{\boldsymbol{\psi}}_h^{\rm eq}}
\newcommand{\BDM}{{\bf BDM}}

\def\avg#1{\left\{\hskip -5pt\left\{#1\right\}\hskip -5pt\right\}}
\def\jump#1{\left[\hskip -4pt\left[#1\right]\hskip -5pt\right]}

\begin{document}
\date{\today}
\maketitle
\begin{abstract}
In this article, goal-oriented a posteriori error estimation for the biharmonic plate bending problem is considered. The error for approximation of goal functional is represented by an estimator which combines dual-weighted residual method and equilibrated moment tensor.
An abstract unified framework for the goal-oriented a posteriori error estimation is derived. In particular, $C^0$ interior penalty and discontinuous Galerkin finite element methods are employed for practical realization. The abstract estimation is based on 
equilibrated moment tensor and potential reconstruction that provides a guaranteed upper bound for the goal error. Numerical experiments are performed to illustrate the effectivity of the estimators.
\end{abstract}


\section{Introduction}\label{intro_Goal_bih}
Adjoint-based goal-oriented a posteriori error estimation is an efficient tool for numerical approximation of many engineering problems since they provide relevant information about error in a quantity of interest rather than error estimates derived in some norm or semi-norm.
The goal-oriented a posteriori error estimation was initially proposed by Becker and Rannacher \cite{Beck_Ran_opt_a_post_FE_01} and by Prudhomme and Oden \cite{Prud_Oden_goal_pointwise_19,Oden_Prudh_goal_or_FE_01} using dual-weighted residual (DWR) method, see \cite{Gil_Sul_goal_02,Bangerth_Ranacher_book_03,Burg_Naza_15_Goal_adap,Mom_Stev_09_Goal_adap_conv} for subsequent works. Some of the popular approaches on goal-oriented a posteriori error estimation are multi-objective goal functional error estimation of \cite{Brum_Zhuk_Zwie_17_goal_multi_obj,Endt_Wick_17_goal_partition_of_unity,Hart_08_Multi_target}, the constitutive relation error (CRE) of \cite{Lad_Stric_upper_bd_08,Lad_Pled_Cham_goal_13,Rey_Rey_Gos_a_post_DD_14,Rey_Rey_Gos_a_post_DD_GO_15,Rey_Rey_Gos_a_post_DD_16}, enhanced least-squares finite element methods of \cite{Chau_Cyr_Liu_14_goal_least_sq}, combination of DWR and equilibrated flux of \cite{Mozo_Prudh_15}, guaranteed bounds based on equilibrated flux of \cite{GM_MV_SY_19,Ains_Ran_12_goal}.

A traditional a posteriori error analysis hinges on the computation of residual \cite{Bangerth_Ranacher_book_03}
\begin{equation}
	a(u-u_h,v)=l(v)-a(u_h,v)=:\rho(u_h)(v)
\end{equation}
for some bilinear form $a(\bullet,\bullet)$ and linear form $l$ associated to an elliptic partial differential equations (PDEs) with $u$ and $u_h$ being its weak and Galerkin solutions, and $v$ being a test function. Incorporating goal-functional $Q$ with a dual problem $a(v,z)=Q(v)$ for all test functions $v$, and using Galerkin orthogonality, we have
\begin{equation}
	Q(u-u_h)=a(u-u_h,z)=a(u-u_h,z-v_h)=\rho(u_h)(z-v_h)    
\end{equation}
for all discrete functions $v_h$.
This approach involves computation of estimator weighted with the solution related to dual problem. The residual based estimator can be chosen for primal problem, and for computational purpose the solution $z$ of dual problem can be chosen as solution obtained in some finer discretization space. However, most of the cases the estimators involve unknown constants, hence they do not provide guaranteed a posteriori estimator. To obtain a guaranteed a posteriori estimator often one incorporates equilibrated flux for second-order PDEs, see \cite{GM_MV_SY_19,Lad_Pled_Cham_goal_13}. Many research has been done in the direction of goal-oriented a posteriori estimation for second-order PDEs. However, to my knowledge, there are hardly very few results (except \cite{Gon_Gom_Dev_14_Goal_Bih_hpDG}) on goal-oriented a posteriori error estimation for fourth-order PDEs. A $hp$-discontinuous Galerkin DWR based goal error estimation has been proposed by \cite{Gon_Gom_Dev_14_Goal_Bih_hpDG} for biharmonic problem and applied to describe the displacement of a thin and isotropic homogeneous plate and the stream function formulation of the Stokes fluid problem describing the flow of a viscous fluid around a flat plate. 

The main purpose of the article is to develop a unified framework for goal-oriented a posteriori error estimation for a model linear biharmonic problem. We consider goal functional of the form $Q(u)=(\tf,u)$ for some weight function $\tf\in\lt$. In practical applications, this can be applied to approximate the goal functional governed by mean deflection around a specified zone and point deflection at some point (in regularized form). In this article, the present framework is applied to (but not limited to) $C^0$ interior penalty and discontinuous Galerkin finite element methods. We also establish a goal estimator which combines the DWR method and equilibrated moment tensor for primal and dual problems. Finally, a unified guaranteed a posteriori error estimation is derived using potential reconstruction and equilibrated moment tensor that is significantly different from the DWR method of \cite{Gon_Gom_Dev_14_Goal_Bih_hpDG}.

The organization of the paper is as follows: Section~\ref{sec:notation} is devoted to notation and preliminaries. Section~\ref{sec:problem} introduces model problem and some useful results. Section~\ref{sec:goal_est} establishes a posteriori error estimates for the goal functional in an abstract framework. Section~\ref{sec:discretization} then considers finite element discretization for the approximation of solution and address some applications of the abstract framework. Finally, in Section~\ref{sec:numer} some numerical experiments are performed to substantiate the theoretical results.

\section{Setting}\label{sec:notation}
Throughout the paper, standard notation on Lebesgue and Sobolev spaces and their norms are employed.
The standard semi-norm and norm on $H^{s}(\Omega)$ (resp. $W^{s,p} (\Omega)$) for $s>0$ are denoted by $|\bullet|_{s}$ and $\|\bullet\|_{s}$ (resp. $|\bullet|_{s,p}$ and $\|\bullet\|_{s,p}$ ). We refer $H^{-m}(\Omega)$ to be the dual space of $H^{m}_0(\Omega)$ with $\langle\bullet,\bullet\rangle_{m}$ denoting the duality product and if $m=2$ we often denote the duality product simply by $\langle\bullet,\bullet\rangle$. Further, let $\boldsymbol{H}(\Div,\Omega)$ be the Hilbert space of vector fields $\underline{\boldsymbol{q}}\in[\lt]^2$ such that $\nabla{\cdot} \underline{\boldsymbol{q}}\in\lt$. Matrix valued functions in $[\lt]^{2\times 2}$ is denoted by $\underline{\underline{\boldsymbol{q}}}=(q_{ij})_{i,j=1}^2$ and the inner-product is $(\underline{\underline{\boldsymbol{p}}},\underline{\underline{\boldsymbol{q}}})=\int_{\Omega} \underline{\underline{\boldsymbol{p}}}:\underline{\underline{\boldsymbol{q}}}\dx$, where $\underline{\underline{\boldsymbol{p}}}:\underline{\underline{\boldsymbol{q}}}=\sum_{i,j=1}^2p_{ij}q_{ij}$. Moreover, we introduce the Hilbert space
\begin{align*}
	\underline{\underline{\boldsymbol{H}}}:=\left\{\underline{\underline{\boldsymbol{q}}}\in \boldsymbol{H}(\Div,\Omega)^2 : \nabla{\cdot}\underline{\underline{\boldsymbol{q}}}\in\boldsymbol{H}(\Div,\Omega)\right\}.
\end{align*}
Finally, we refer to $D^2v:=(\partial^2 v/\partial x_i\partial x_j)_{i,j=1}^2$ as the matrix of second order partial derivatives of a function $v\in H^2(\Omega)$. The set of all symmetric $2\times 2$ matrix valued functions are denoted by $\symmat$.

Let  $\cT_h$ be a shape-regular \cite{Braess} triangulation of the bounded polygonal Lipschitz domain $\Omega\subset\mR^2$ into closed triangles. 
The set of all internal vertices (resp. boundary vertices) and  interior edges (resp.  boundary edges)  of the triangulation $\cT_h$ are denoted by $\cN_h (\Omega)$ (resp.  $\cN_h(\partial\Omega)$) and $\cE_h (\Omega)$ (resp. $\cE_h (\partial\Omega)$).
Define a piecewise constant mesh function $h_{\cT_h}(x)=h_K={\rm diam} (K)$ for all $x \in K$, $ K\in \cT_h$, and set $h:=\max_{K\in \cT_h}h_K$. Also define a piecewise constant edge-function on $\cE_h:=\cE_h(\Omega)\cup \cE_h(\partial\Omega)$ by $h_{\cE_h}|_{\e}=h_{\e}={\rm diam}(e)$ for any $e\in \cE_h$. Set of all edges of $K$ is denoted by $\cE_h(K)$. Note that for a shape-regular family, there exists a positive constant $C$ independent of $h$ such that any $K\in\cT_h$ and any $e\in \cE_h(K)$ satisfy $Ch_K\leq h_{\e}\leq h_K$. Let $\mP_k( K)$ denote the set of all polynomials of degree less than or equal to $k$ and 
\[\mP_k(\cT_h):=\left\{\varphi\in L^2(\Omega):\,\forall K\in\cT_h,\varphi|_{K}\in \mP_k(K)\right\}.\]
The $L^2(\Omega)$ projection onto $\mP_k(\cT_h)$ is denoted by $\Pi_k$.
For a nonnegative integer $s$, define the broken Sobolev space for the subdivision $\cT_h$  as
\begin{equation*}
	H^s(\cT_h)=\left\{\varphi\in\lt: \varphi|_K\in H^{s}(K)\fl K\in \cT_h \right\},
\end{equation*}
with the broken Sobolev semi-norm  $|\bullet|_{H^s(\cT_h)}$ and norm $\| \bullet\|_{H^s(\cT_h)}$ defined by
\begin{equation*}
	|\varphi|_{H^s(\cT_h)}=\bigg{(}\sum_{K\in\cT_h} |\varphi|_{H^{s}(K)}^2\bigg{)}^{1/2}\text{ and }
	\|\varphi\|_{H^s(\cT_h)}=\bigg{(}\sum_{K\in\cT_h}\|\varphi\|_{H^{s}(K)}^2\bigg{)}^{1/2}.
\end{equation*}
Define the jump $\jump{\varphi}_{\e}=\varphi|_{K_+}-\varphi|_{K_-}$ and the average $\left\{\hskip -4pt\left\{\varphi\right\}\hskip -4pt\right\}_{\e}=\half\left(\varphi|_{K_+}+\varphi|_{K_-}\right)$ across the interior edge $e$ of $\varphi\in H^1(\cT_h)$ of the adjacent triangles  $K_+$ and $K_-$. Extend the definition of the jump and the average to an edge lying on boundary by $\jump{\varphi}_{\e}=\varphi|_{\e}$ and $\left\{\hskip -4pt\left\{\varphi\right\}\hskip -4pt\right\}_{\e}=\varphi|_{\e}$ { when $e$ belongs to the set of boundary edges $\cE_h(\partial\Omega)$}.
For any vector function, jump and average are understood componentwise.

There exist real numbers $\Ctr$ and $\Ctrc$ independent of $h$ such that the following discrete and continuous trace inequalities hold for all $K\in\cT_h$ and $e\in\cE_h$ (see \cite[Lemma~1.46 and 1.49]{Piet_Ern_12_DG_book})
\begin{align}
	\|v\|_{\e}&\leq \Ctr h_{\e}^{-1/2} \|v\|_K\quad\forall v\in \mP_k(K),\label{dis_trace}\\ 
	\|v\|_{\partial K}&\leq C_{\text{tr},\text{c}}(h_K^{-1}\|v\|_T^2+ h_K\|\nabla v\|_K^2)^{1/2}\quad\forall v\in H^1(K).
\end{align}

The positive constants $C$ appearing in the inequalities denote generic constants which do not depend on the mesh-size. The notation $a\lesssim b$ means that there exists a generic constant $C$ independent of the mesh parameters such that $a \leq Cb$.

\section{Model problem}\label{sec:problem}
In this article, we are interested in general linear fourth-order boundary-value problems, but the results can be extended to more general situations. For the sake of simplicity of presentation, we restrict ourselves to a simple model problem. Consider the biharmonic equation with clamped boundary conditions
\begin{subequations}\label{eqn:bih}
	\begin{align}
		&\Delta^2 u=f \quad\text{ in }\Omega,\\
		&u=0=\frac{\partial u}{\partial\nu}\quad \text{ on } \partial\Omega,
	\end{align}
\end{subequations}
where $\Delta^2 u=\Delta (\Delta u)$ and the source term $f$. Define $V:=\hto$.
The weak formulation is given by: for  $f\in H^{-2}(\Omega)$, find $u\in V$ such that
\begin{equation}\label{weak_prim}
	(D^2u,D^2v)=(f,v)\fl v\in V.
\end{equation}

In this article, we are interested on the following goal functional
\begin{equation}\label{goal_func}
	Q(u)=(\tf,u)
\end{equation}
for some chosen weight function $\tf\in L^2(\Omega)$. We analyze the above goal functional by a dual problem of \eqref{eqn:bih} that consists in finding $\tu:\Omega\to\mR$ such that
\begin{subequations}\label{eqn:bih_dual}
	\begin{align}
		&\Delta^2 \tu=\tf \quad\text{ in }\Omega,\\
		&\tu=0=\frac{\partial \tu}{\partial\nu}\quad \text{ on } \partial\Omega,
	\end{align}
\end{subequations}
and the weak formulation seeks $\tu\in V$ such that
\begin{equation}\label{weak_dual}
	(D^2\tu,D^2v)=(\tf,v)\fl v\in V.
\end{equation}
Existence and uniqueness of the weak solution of both the primal and dual problems \eqref{weak_prim} and \eqref{weak_dual} follow
from Riesz representation theorem.

We state two definitions which are essential for establishing a posteriori error estimation.
\begin{definition}[Potential reconstruction]\label{Defn_Potn_pecons}
	We call a potential reconstruction any function $s_h$ (resp. $\ts_h$) constructed from $u_h$ (resp. $\tu_h$) which satisfies 
	\begin{align}
		s_h\in \hto\cap C^1(\bar{\Omega})\quad (\text{resp. } \ts_h\in \hto\cap C^1(\bar{\Omega})).\label{poten_cont}
	\end{align}	
\end{definition}

Here and throughout the paper we consider the $\Div\Div$ operator in distributional sense, i.e., for $\underline{\underline{\tau}}\in \symmat$
\begin{equation}\label{eqn:defn_DivDiv}
	\langle \Div\Div \underline{\underline{\tau}},w\rangle=\int_\Omega \underline{\underline{\tau}} : D^2 w\dx, \quad \fl w\in\hto.
\end{equation}

\begin{definition}[Equilibrated moment tensor]\label{Defn_moment_tensor} Let $f_h\in H^{-2}(\Omega)$ (resp. $\tf_h\in H^{-2}(\Omega)$ ).  Any matrix valued function $\eqflux\in \symmat$ (resp. $\teqflux\in \symmat$) which satisfies
	
	\begin{subequations}
		\begin{align}
			&\langle\Div\Div \eqflux,w\rangle=\langle f_h,w\rangle \fl w\in\hto \label{recons_flux}\\
			(\text{resp. }  &\langle\Div\Div \teqflux,w\rangle=\langle\tf_h,w\rangle \fl w\in\hto)
		\end{align}
	\end{subequations}
	is called an equilibrated moment tensor.
\end{definition}	


The following Prager--Synge type energy principle hold: 
\begin{lemma}[Two-energies principle for biharmonic equation]\label{Prager_Synge}
	Let $f_h\in H^{-2}(\Omega)$ and $\hat{u}\in\hto$ be the solution of the biharmonic equation
	\begin{equation}\label{bih_fh}
		(D^2 \hat{u},D^2 w)=\langle f_h,w\rangle\fl w\in\hto.
	\end{equation}	
	For $v\in\hto$, the tensor $\eqflux\in\symmat$ defined in Definition~\ref{Defn_moment_tensor} satisfies \cite{Prag_Syng_47,Bra_Pech_Sch_20_equi_flux_bih}
	\begin{equation}\label{eqn_prag_syn}
		\|D^2(\hat{u}-v)\|^2+\|D^2\hat{u}-\eqflux\|^2=\|D^2 v-\eqflux\|^2.
	\end{equation}
	Moreover, let  $u\in\hto$ (resp. $\tu\in\hto$) be the solution of \eqref{weak_prim} (resp. \eqref{weak_dual}). Then for any $v\in\hto$ (resp. $\tilde{v}\in\hto$), the following also holds
	\begin{align}
		&\|D^2(u-v)\|^2+\|D^2u-\eqflux\|^2=\|D^2 v-\eqflux\|^2+2\langle f-f_h,u-v\rangle.\label{eqnPrim_prag_syn_osc}\\
		(\text{resp. } &\|D^2(\tu-\tilde{v})\|^2+\|D^2\tu-\teqflux\|^2=\|D^2 \tilde{v}-\teqflux\|^2+2\langle\tf-\tf_h,\tu-\tilde{v}\rangle.)\label{eqnDual_prag_syn_osc}
	\end{align}
\end{lemma}	
\begin{proof}
	For $v\in\hto$, adding and subtracting $\hat{u}$ then expanding, we have
	\begin{align}\label{Prag_Syn_id}
		\|D^2 v-\eqflux\|^2=\|D^2(v-\hat{u})\|^2+\|D^2\hat u-\eqflux\|^2+2(D^2(v-\hat{u}),D^2\hat{u}-\eqflux).
	\end{align}
	We use \eqref{bih_fh}  and the definition of $\eqflux$ in \eqref{recons_flux} to obtain
	\begin{equation}
		(D^2(v-\hat{u}),D^2\hat{u}-\eqflux)=(D^2(v-\hat{u}),D^2\hat{u})-(D^2(v-\hat{u}),\eqflux)=\langle f_h,v-\hat{u}\rangle-\langle f_h,v-\hat{u}\rangle=0.
	\end{equation} 
	The above two equations yield \eqref{eqn_prag_syn}.
	
	The proof of \eqref{eqnPrim_prag_syn_osc} follows from the identity \eqref{Prag_Syn_id} (replacing $\hat{u}$ by $u$) with \eqref{weak_prim} and \eqref{recons_flux}
	\begin{align*}
		(D^2(v-u),D^2u -\eqflux)&=(D^2(v-u),D^2u)-(D^2(v-u),\eqflux)=\langle f-f_h,v-u\rangle.
	\end{align*}
	The expression~\eqref{eqnDual_prag_syn_osc} is proved similarly by exploring dual problem.
\end{proof}

\section{Goal-oriented error estimate}\label{sec:goal_est}

Choosing $v=u$ in \eqref{weak_dual} and $v=\tu$ in \eqref{weak_prim}, the following primal-dual equivalence relation holds
\begin{equation}\label{prim_dual_eq_pel}
	Q(u)=(\tf,u)=(D^2\tu,D^2 u)=(D^2u,D^2 \tu)=\langle f,\tu\rangle.
\end{equation}

In the following subsections, goal functional is approximated and some error representations are presented.

\subsection{Some residual type goal error estimations}
In this subsection, goal error is represented by an estimator and by a remainder term. For any edge $e\in \cE_h$, the outward unit normal across the edge is denoted by $\bfn$ and unit tangent along the edge is denote by $\boldsymbol{\tau}_e$. Define $\partial_n v:=\nabla v{\cdot}\bfn, D^2_{nn}v:=\bfn^TD^2v\bfn$ and $\eqfnn:=\bfn^T\eqflux\bfn, \eqfnt:=\boldsymbol{\tau}_e^T\eqflux\bfn$.

\begin{theorem}[Error representation of goal functional]\label{thm:res_est}
	Let $u$ and $\tu$ respectively be the solutions of \eqref{eqn:bih} and \eqref{eqn:bih_dual}. Let $u_h$ and $\tu_h\in\mP_k(\cT_h)$ be arbitrary piecewise polynomial functions. Let $s_h$ and $\ts_h$ be the potential reconstructions of Definition~\ref{Defn_Potn_pecons}, and $\eqflux$ and $\teqflux$ be the equilibrated moment tensors of Definition~\ref{Defn_moment_tensor} constructed from $u_h$ and $\tu_h$ respectively. Then the goal error is expressed as
	\begin{align}\label{goal_err_pes}
		Q(u)-Q(u_h)=\eta_h(u_h,\tu_h;\eqflux,\teqflux)+\cR_h(u,\tu,f;u_h,\tu_h),
	\end{align}
	where the estimator is given by
	\begin{align}
		\eta_h:&=\langle f,\tsh\rangle -\sik \eqflux:D^2\tsh\dx+ \sik(\eqflux-D^2u_h):\teqflux\dx \notag\\
		&\quad+\sum_{\e\in\cE_h}\int_{\e} \jump{\partial_{\tau} u_h}_{\e}\teqfnt\dss+ \sum_{\e\in\cE_h}\int_e\jump{\partial_n u_h}_{\e}\teqfnn\dss+\sum_{\e\in\cE_h}\int_{\e} \jump{u_h}_{\e}\Div\teqflux{\cdot}\bfn\dss,\label{goal_main_est}
	\end{align}
	with remainder term
	\begin{align}
		\cR_h:&=\langle f-\Div\Div\eqflux ,\tu-\tsh\rangle + \sum_{K\in\cT_h}\int_K(\eqflux-D^2u_h):(D^2\tu-\teqflux)\dx\notag\\
		&\quad+ \sum_{\e\in\cE_h}\int_e\jump{\partial_{\tau} u_h}_{\e}(D^2_{n\tau}\tu-\teqfnt)\dss+\sum_{\e\in\cE_h}\int_e\jump{\partial_n u_h}_{\e}(D^2_{nn}\tu-\teqfnn)\dss\notag\\
		&\qquad+\sum_{\e\in\cE_h}\int_{\e} \jump{u_h}_{\e} (\Div(D^2\tu)-\Div\teqflux){\cdot}\bfn\dss.\label{goal_rem_est}
	\end{align}
\end{theorem}

\begin{proof}
	The primal-dual equivalence relation~\eqref{prim_dual_eq_pel} and the definition of goal functional  lead to the goal error representation
	\begin{align}\label{eqn:goal_rep_ftf}
		Q(u)-Q(u_h)=\langle f,\tu\rangle-(\tf,u_h).
	\end{align}	
	The dual problem~\eqref{eqn:bih_dual} with regularity of $\tu$ and successive application of integration by parts yield for the above second term as 
	\begin{align}
		(\tf,u_h)&=\sik u_h\Delta^2\tu \dx=\sik u_h\Div\Div(D^2\tu)\dx\notag\\
		&=-\sik\nabla u_h{\cdot}\Div(D^2\tu)\dx+\sum_{K\in\cT_h}\int_{\partial K} u_h\Div(D^2\tu){\cdot}\bn\dss\notag\\
		&=\sik D^2u_h:D^2\tu\dx-\sum_{K\in\cT_h}\int_{\partial K} \nabla u_h{\cdot} D^2\tu\bn\dss+\sum_{K\in\cT_h}\int_{\partial K} u_h\Div(D^2\tu){\cdot}\bn\dss.
	\end{align}
	Expressing gradient in tangent-normal direction as $\nabla u_h=\partial_{\tau} u_h\boldsymbol{\tau}_e+\partial_n u_h\bfn$, we have from the above equation
	\begin{align}
		(\tf,u_h)&=\sik D^2u_h:D^2\tu\dx-\sum_{\e\in\cE_h}\int_{\e} \jump{\partial_{\tau} u_h}_{\e} D_{n\tau}\tu\dss-\sum_{\e\in\cE_h}\int_{\e} \jump{\partial_n u_h}_{\e} D_{nn}\tu\dss\notag\\
		&\qquad+\sum_{\e\in\cE_h}\int_{\e} \jump{u_h}_{\e}\Div(D^2\tu){\cdot}\bfn\dss.\label{eqn:tf_rep}
	\end{align}
	The above two displayed equations \eqref{eqn:goal_rep_ftf} and \eqref{eqn:tf_rep} lead to 
	\begin{align}\label{simple_err_pep}
		Q(u)-Q(u_h)&=\langle f,\tu\rangle-\sik D^2u_h:D^2\tu\dx+\sum_{\e\in\cE_h}\int_{\e} \jump{\partial_{\tau} u_h}_{\e} D_{n\tau}\tu\dss+\sum_{\e\in\cE_h}\int_{\e} \jump{\partial_n u_h}_{\e} D_{nn}\tu\dss\notag\\
		&\qquad\qquad-\sum_{\e\in\cE_h}\int_{\e} \jump{u_h}_{\e}\Div(D^2\tu){\cdot}\bfn\dss..
	\end{align}
	Introducing the equilibrated moment tensor $\eqflux$ and $\teqflux$ of Definition~\ref{Defn_moment_tensor} in the first two terms of the above equation~\eqref{simple_err_pep} yields
	\begin{align}
		&\langle f,\tu\rangle-\sik D^2u_h:D^2\tu\dx\notag\\
		&=\langle f-\Div\Div\eqflux,\tu\rangle+\int_{\Omega}\eqflux:D^2\tu\dx-\sik D^2u_h:D^2\tu\dx\notag\\
		&=\langle f-\Div\Div\eqflux,\tsh\rangle+\langle f-\Div\Div\eqflux,\tu-\tsh\rangle\notag\\
		&\qquad+\sik(\eqflux-D^2u_h):\teqflux\dx+\sik(\eqflux-D^2u_h):(D^2\tu-\teqflux)\dx.\label{est_term}
	\end{align}
	First two terms in the above equation can be written as
	\begin{align}
		\langle f-\Div\Div\eqflux,\tsh\rangle+\langle f-\Div\Div\eqflux,\tu-\tsh\rangle=\langle f,\tsh\rangle -(\eqflux,D^2\tsh)+\langle f-f_h,\tu-\tsh\rangle.
	\end{align}
	Introducing equilibrated moment tensors of tangent-normal directions in the third  and fourth terms of \eqref{simple_err_pep} lead to
	\begin{align}
		&\sum_{\e\in\cE_h}\int_{\e} \jump{\partial_{\tau} u_h}_{\e} D_{n\tau}\tu\dss+\sum_{\e\in\cE_h}\int_{\e} \jump{\partial_n u_h}_{\e} D_{nn}\tu\dss=\sum_{\e\in\cE_h}\int_{\e} \jump{\partial_{\tau} u_h}_{\e}\teqfnt\dss+\sum_{\e\in\cE_h}\int_{\e} \jump{\partial_n u_h}_{\e}\teqfnn\dss\notag\\
		&\qquad+\sum_{\e\in\cE_h}\int_{\e} \jump{\partial_{\tau} u_h}_{\e} (D_{n\tau}\tu-\teqfnt)\dss+\sum_{\e\in\cE_h}\int_{\e} \jump{\partial_n u_h}_{\e} (D_{nn}\tu-\teqfnn)\dss.
	\end{align}
	Introducing equilibrated moment tensor of normal direction in the last term of \eqref{simple_err_pep} leads to
	\begin{align}\label{jump_div}
		\sum_{\e\in\cE_h}\int_{\e} \jump{u_h}_{\e}\Div(D^2\tu){\cdot}\bfn\dss=\sum_{\e\in\cE_h}\int_{\e} \jump{u_h}_{\e}\Div\teqflux{\cdot}\bfn\dss+\sum_{\e\in\cE_h}\int_{\e} \jump{\partial_n u_h}_{\e} (\Div(D^2\tu)-\Div\teqflux){\cdot}\bfn\dss.
	\end{align}	
	The last five displayed equations \eqref{simple_err_pep}-\eqref{jump_div} represent the goal error equation \eqref{goal_err_pes} with the estimator term $\eta_h$ and remainder term $\cR_h$.
\end{proof}


\begin{remark}[Goal estimator 2]\label{rem:goal_est2}
	If $f$ and $\Div\Div \eqflux$ belong to $L^2(\Omega)$, we can replace $s_h$ by $u_h$ in \eqref{est_term}and obtain a simplified estimator of Theorem~\ref{thm:res_est} as
	\begin{align}
		\eta_h:&= (f-\Div\Div \eqflux,\tu_h)+ \sik(\eqflux-D^2u_h):\teqflux\dx\notag\\
		&\quad+\sum_{\e\in\cE_h}\int_{\e} \jump{\partial_{\tau} u_h}_{\e}\teqfnt\dss+ \sum_{\e\in\cE_h}\int_e\jump{\partial_n u_h}_{\e}\teqfnn\dss+\sum_{\e\in\cE_h}\int_{\e} \jump{u_h}_{\e}\Div\teqflux{\cdot}\bfn\dss.\label{simple_goal_est}
	\end{align}
	with remainder term
	\begin{align}
		\cR_h:&=(f-\Div\Div \eqflux,\tu-\tu_h) + \sum_{K\in\cT_h}\int_K(\eqflux-D^2u_h):(D^2\tu-\teqflux)\dx\notag\\
		&\quad+\sum_{\e\in\cE_h}\int_e\jump{\partial_{\tau} u_h}_{\e}(D^2_{n\tau}\tu-\teqfnt)\dss+\sum_{\e\in\cE_h}\int_e\jump{\partial_n u_h}_{\e}(D^2_{nn}\tu-\teqfnn)\dss\notag\\
		&\qquad+\sum_{\e\in\cE_h}\int_{\e} \jump{u_h}_{\e} (\Div(D^2\tu)-\Div\teqflux){\cdot}\bfn\dss.
	\end{align}
	Moreover, the estimator $\eta_h$ of \eqref{simple_goal_est} can be expressed as the sum of local element error contributions:
	\begin{align}
		\eta_h=\sum_{K\in\cT_h}\eta_K=\sum_{K\in\cT_h}(\eta_{{\rm est},K}+\eta_{{\rm jump},K}+\eta_{\mathcal{O},K}),
	\end{align}
	where the local contributions are given by
	\begin{align*}
		\eta_{{\rm est},K}&:=\int_K(\eqflux-D^2u_h):\teqflux\dx,\\
		\eta_{{\rm jump},K}&:=\sum_{\e\in\cE_K}\int_e\gamma_\e\left(\jump{\partial_n u_h}_{\e}\teqfnn+\jump{\partial_{\tau} u_h}_{\e}\teqfnt+\jump{u_h}_{\e}\Div\teqflux{\cdot}\bfn\right)\dss,\\
		\eta_{\mathcal{O},K}&:=\int_K (f-\Div\Div \eqflux)\tuh\dx,
	\end{align*}
	with the indicator function $\gamma_\e=1/2$ for interior edge $\e\in\cE_h({\Omega})$ and $\gamma_\e=1$ for boundary edge $\e\in\cE_h({\partial\Omega})$.
\end{remark}

\subsection{Guaranteed a posteriori error estimate}

In this subsection, we present a guaranteed a posteriori error estimator for the goal-error based on  equilibrated moment tensor and potential reconstruction. An abstract a posteriori estimator is proposed, and is computed later for some finite element methods. Here and throughout this subsection, for given $\eqflux$ and $\teqflux$ belong to $\symmat$, we define $f_h:=\Div\Div \eqflux$ and $\tf_h:=\Div\Div \teqflux$. We proceed first by writing goal-error as:
\begin{lemma}[Goal error equation]
	Let $u$ and $\tu$ respectively be the solution of \eqref{eqn:bih} and \eqref{eqn:bih_dual}. Let $u_h$ and $\tu_h\in\mP_k(\cT_h)$ respectively be arbitrary piecewise polynomial functions. Let $s_h$ and $\ts_h$ be the potential reconstructions of Definition~\ref{Defn_Potn_pecons}, and $\eqflux$  be the equilibrated moment tensors of Definition~\ref{Defn_moment_tensor}. There holds
	\begin{equation}\label{eq_goal_err}
		Q(u)-Q(u_h)=\langle f-f_h,\tu\rangle+(\eqflux-D^2 s_h,D^2\tu)+Q(s_h-u_h).
	\end{equation}
\end{lemma}
\begin{proof}
	From the primal dual equivalence relation \eqref{prim_dual_eq_pel} and Definition~\ref{Defn_moment_tensor}, we obtain
	\begin{align}
		Q(u)&=\langle f,\tu\rangle=\langle f-f_h,\tu\rangle+\langle f_h,\tu\rangle\notag\\
		&=\langle f-f_h,\tu\rangle+\langle\Div\Div \eqflux,\tu\rangle=\langle f-f_h,\tu\rangle+(\eqflux,D^2\tu).
	\end{align}
	Since $s_h\in\hto$, from the weak formulation of dual problem \eqref{weak_dual} with $v=s_h$, we obtain
	\begin{align}
		Q(u_h)=Q(s_h)+Q(u_h-s_h)=(D^2\tu,D^2 s_h)+Q(u_h-s_h)=(D^2 s_h,D^2 \tu)+Q(u_h-s_h).
	\end{align}
	From the above two displayed equations we have
	\begin{equation}
		Q(u)-Q(u_h)=\langle f-f_h,\tu\rangle+(\eqflux-D^2 s_h,D^2\tu)-Q(u_h-s_h).
	\end{equation}
	This completes the proof.
\end{proof}

We apply the principle of classical bounding technique of Ladev{\`e}ze {\it et al.} \cite{Lad_Cham_goal_10,Lad_Stric_upper_bd_08} related to the goal-oriented a posteriori error estimate of elasticity problem. Let 
\begin{equation}\label{avg_flux_pecons}
	\tbfsigh^{\rm m}:=\frac{1}{2}\left(\teqflux+D^2\ts_h^i\right)
\end{equation} 
be the average of moment tensor  of Definition~\ref{Defn_moment_tensor} and Hessian of potential reconstruction of Definition~\ref{Defn_Potn_pecons} corresponding to the dual problem. We denote the following  oscillation terms  by
\begin{equation}\label{defn_data_osc}
	osc^2_{\rm prim}(f,\tu):=|\langle f-f_h,\tu-\ts_h\rangle| \text{ and } osc^2_{\rm dual}(\tf,\tu):=|\langle\tf-\tf_h,\tu-\ts_h\rangle|.
\end{equation}

\begin{theorem}[Abstract goal-oriented a posteriori estimator]
	\label{thm:abs_goal_est} Let $u$ and $\tu$ respectively be the solution of \eqref{eqn:bih} and \eqref{eqn:bih_dual}. Let $u_h$ and $\tu_h\in\mP_k(\cT_h)$ respectively be arbitrary piecewise polynomial functions. Let $s_h$ and $\ts_h$ be the potential reconstructions of Definition~\ref{Defn_Potn_pecons}, and $\eqflux$ and $\teqflux$ be the equilibrated moment tensors of Definition~\ref{Defn_moment_tensor} with $\tbfsigh^{\rm m}$ being the average moment tensor of \eqref{avg_flux_pecons}. There holds
	\begin{align}
		&\left|Q(u)-Q(u_h)-\left(\eqflux-D^2 s_h,\tbfsigh^{\rm m}\right)\right|\notag\\
		&\leq  \|D^2 s_h-\eqflux\|\left(\half\|D^2 \ts_h-\teqflux\|+osc_{\rm dual}(\tf,\tu)\right)+|\langle f-f_h,\ts_h\rangle+Q(s_h-u_h)|+osc^2_{\rm prim}(f,\tu).\label{abs_goal_est}
	\end{align}
\end{theorem}

\begin{proof}
	Adding and subtracting the average moment tensor $\tbfsigh^{\rm m}$ in \eqref{eq_goal_err}, we obtain
	\begin{equation}\label{goal_err_mean}
		Q(u)-Q(u_h)-\left(\eqflux-D^2 s_h,\tbfsigh^{\rm m}\right)=\langle f-f_h,\tu\rangle+(\eqflux-D^2 s_h,D^2\tu-\tbfsigh^{\rm m})+Q(s_h-u_h).
	\end{equation}
	From the definition of \eqref{avg_flux_pecons} and expanding, we have
	\begin{align*}
		\|D^2\tu-\tbfsigh^{\rm m}\|^2=\frac{1}{4}\|D^2(\tu-\ts_h)\|^2+\frac{1}{4}\|D^2 \tu-\teqflux\|^2+\half(D^2(\tu-\ts_h),D^2 \tu-\teqflux).
	\end{align*}
	Apply integration by parts twice to obtain
	\begin{equation*}
		(D^2(\tu-\ts_h),D^2 \tu-\teqflux)=\langle\tu-\ts_h,\tf-\tf_h\rangle.
	\end{equation*}
	The above two equations and \eqref{eqnDual_prag_syn_osc} with $\tilde{v}=\ts_h$ imply
	\begin{align}\label{dual_est_bdd}
		\|D^2\tu-\tbfsigh^{\rm m}\|^2=\frac{1}{4}\|D^2 \ts_h-\teqflux\|^2+\langle\tu-\ts_h,\tf-\tf_h\rangle\leq \left(\half\|D^2 \ts_h-\teqflux\|+osc_{\rm dual}(\tf,\tu)\right)^2.
	\end{align}
	Apply Schwarz inequality in the right hand side of \eqref{goal_err_mean} and use \eqref{dual_est_bdd} to obtain
	\begin{align*}
		&|Q(u)-Q(u_h)-\left(\eqflux-D^2 s_h,\tbfsigh^{\rm m}\right)|\\
		&\leq  \|D^2 s_h-\eqflux\|\left(\half\|D^2 \ts_h-\teqflux\|+osc_{\rm dual}(\tf,\tu)\right)+|\langle f-f_h,\tu\rangle+Q(s_h-u_h)|\\
		&\leq  \|D^2 s_h-\eqflux\|\left(\half\|D^2 \ts_h-\teqflux\|+osc_{\rm dual}(\tf,\tu)\right)+|\langle f-f_h,\ts_h\rangle+Q(s_h-u_h)|+osc_{\rm prim}(f,\tu).
	\end{align*}
	This completes the proof.
\end{proof}

\begin{remark}
	A use of the the identity \eqref{goal_err_mean}, the above goal a posteriori error estimation \eqref{abs_goal_est} can be rewritten as 
	\begin{align}\label{abs_goal_est_simple}
		\left|Q(u)-Q(s_h)-\left(\eqflux-D^2 s_h,\tbfsigh^{\rm m}\right)\right|&\leq  \|D^2 s_h-\eqflux\|\left(\half\|D^2 \ts_h-\teqflux\|+osc_{\rm dual}(\tf,\tu)\right)\notag\\
		&\quad+|\langle f-f_h,\ts_h\rangle|+osc_{\rm prim}^2(f,\tu).
	\end{align}
\end{remark}

\begin{remark}
	The oscillation terms can be computed as follows. The triangle inequality and Lemma~\ref{Prager_Synge} imply
	\begin{align}\label{nc_est}
		\|\tu-\ts_h\|_2\leq \|\tu-\hat{\tu}\|_{2}+\|\hat{\tu}-\ts_h\|_{2,h}\leq \|\tf_h-\tf\|_{-2}+\|D^2 \ts_h-\teqflux\|.
	\end{align}
	This leads to a bound for data oscillation defined in \eqref{defn_data_osc},
	\begin{align}
		osc^2_{\rm dual}(\tf,\tu)&:=|(\tf-\tf_h,\tu-\ts_h)|\leq \|\tf-\tf_h\|_{-2}\|\tu-\ts_h\|_2 \notag\\
		&\leq \|\tf-\tf_h\|_{-2}\left(\|\tf-\tf_h\|_{-2}+\|D^2 \ts_h-\teqflux\|\right).\label{osc_d}
	\end{align}
	Similarly, the second data oscillation in \eqref{defn_data_osc} can be bounded as
	\begin{align}\label{osc_k}
		osc^2_{\rm prim}(f,\tu)\leq\|f-f_h\|_{-2}\left(\|\tf-\tf_h\|_{-2}+\|D^2 \ts_h-\teqflux\|\right).
	\end{align}
	We observe that if there are no data oscillations for primal and dual problems, then $osc_{\rm prim}(\tf,\tu)=0$ and $osc_{\rm dual}(f,\tu)=0$. Then the abstract a posteriori estimator \eqref{abs_goal_est} yields the simplified form:
	\begin{equation}
		\left|Q(u)-Q(u_h)-\left(\eqflux-D^2 s_h,\tbfsigh^{\rm m}\right)\right|
		\leq  \half\|D^2 s_h-\eqflux\|\|D^2 \ts_h-\teqflux\|+|Q(s_h-u_h)|.
	\end{equation}
\end{remark}



\section{Discretization of the biharmonic equation}\label{sec:discretization}
In this section, two nonstandard finite element methods are discussed in order to realize the estimator found in Section~\ref{sec:goal_est}. First, finite element approximation is introduced, then some procedures are described to obtain potential reconstruction of Definition~\ref{Defn_Potn_pecons} and equilibrated moment tensor of Definition ~\ref{Defn_moment_tensor}. Here and rest of the article, we assume $f,\tf\in\lt$ and $k\geq 2$.
\subsection{$C^0$IPDG method}\label{sec:c0ip}
We obtain approximate solution by $C^0$ interior penalty method ($C^0$IPDG ); see \cite{BGS10, Bra_Pech_Sch_20_equi_flux_bih}. Define the polynomial space for $C^0$IPDG by
\begin{equation*}
	V_h^k:=\{v_{\ip}\in C^0(\Omega)\; | \; v_{\ip}|_K\in \mP_k(K),\; K\in\cT_h\}.
\end{equation*}

Define the bilinear form $a_{\ip}(\cdot,\cdot): V_h^k\times V_h^k\to \mR$ by 
\begin{align}
	a_{\ip}(u_{\ip},v_{\ip})&:=\sum_{K\in\cT_h}\int_K D^2 u_{\ip}:D^2 v_{\ip}\dx-\sum_{\e\in\cE_{h}}\int_\e \left(\jump{\partial_n u_{\ip}}_{\e}\avg{D^2 v_{\ip,nn}}_{\e}+\avg{D^2u_{\ip,nn}}_{\e}\jump{\partial_nv_{\ip}}_{\e}\right)\dss\notag\\
	&\quad+\sum_{\e\in\cE_h}\frac{\sigma}{h_\e}\int_\e \jump{\partial_n u_{\ip}}_{\e}\jump{\partial_nv_{\ip}}_{\e}\dss
\end{align}
where $\sigma$ is  large positive penalty parameter. Define linear forms for primal and dual problems as
\begin{equation}
	l_{\ip}(v_{\ip}):=\sik fv_{\ip}\dx \text{ and } \tilde{l}_{\ip}(v_{\ip}):=\sik \tf v_{\ip}\dx\fl v_{\ip}\in V_h^k.
\end{equation}
The $C^0$IPDG method for \eqref{weak_prim} seeks $u_{\ip}\in V_h^k$ such that
\begin{equation}\label{c0ip_prim}
	a_{\ip}(u_{\ip},v_{\ip})=l_{\ip}(v_{\ip})\fl v_{\ip}\in V_h^k,
\end{equation}
and $C^0$IPDG method for the dual problem \eqref{weak_dual} seeks $\tu_{\ip}\in V_h^k$ such that
\begin{equation}\label{c0ip_dual}
	a_{\ip}(\tu_{\ip},v_{\ip})=\tilde{l}_{\ip}(v_{\ip})\fl v_{\ip}\in V_h^k.
\end{equation}
The discretization error is measured by the mesh-dependent norm
\begin{equation*}
	\|v\|_{\ip}^2:=\sum_{K\in\cT_h}\|D^2 v\|_{0,K}^2+\sum_{\e\in\cE_h}\frac{\sigma}{h_\e}\|\jump{\partial_n v}_{\e}\|_{0,\e}^2 \fl v\in V_h^k+\hto.
\end{equation*}
It is well known that for sufficiently large $\sigma=O((k+1)^2)$, there exists a positive constant $\beta$ such that the following coercivity result holds:
\begin{equation*}
	a_{\ip}(v_{\ip},v_{\ip})\geq \beta\|v_{\ip}\|_{\ip}^2 \fl v_{\ip}\in V_h^k.
\end{equation*}
Also, the bilinear form $a_{\ip}(\cdot,\cdot)$ is continuous, i.e., $|a_{\ip}(v,w)|\leq C\|v\|_{\ip}\|w\|_{\ip}$ for all $ v,w\in V_h^k$. 
The boundedness and coercivity of $a_{\ip}(\cdot,\cdot)$, and continuity of $l_{\ip}$ and $\tilde{l}_{\ip}$ lead to existence and uniqueness of solution of primal and dual problems \eqref{c0ip_prim}-\eqref{c0ip_dual} by Lax-Milgram lemma.

The estimator of Theorem~\ref{thm:res_est} is computed by the construction of equilibrated moment tensors $\eqflux$ and $\teqflux$ of Definition~\ref{Defn_moment_tensor}, and potential reconstructions $s_h$ and $\ts_h$ of Definition~\ref{Defn_Potn_pecons}. Their constructions are outlined below:

\noindent\textbf{Construction of equilibrated moment tensor.}
We follow \cite{Bra_Pech_Sch_20_equi_flux_bih} for the construction of equilibrated moment tensor. Define the symmetric piecewise polynomial tensor fields of order $k-1$ with continuous normal-normal component $\underline{\underline{\tau}}_{h,nn}=\bfn^T\tauh \bfn$ by
\begin{align}
	\Mh:=\{\tauh\in [L^2(\Omega)]^{2\times2}_{\rm sym}\; |\; &\tauh\in[\mP_{k-1}(K)]^{2\times 2}_{\rm sym}, K\in\cT_h,\notag\\
	& \underline{\underline{\tau}}_{h,nn} \text{ is continuous at interelement boundaries} \}.
\end{align}
Each $\tauh\in \Mh$ is uniquely defined by the degrees of freedom (see \cite{Comodi_89_HHJ,Bra_Pech_Sch_20_equi_flux_bih})
\begin{align}\label{EquTensorDof}
	&\int_\e\underline{\underline{\tau}}_{h,nn}q_\e\dss, \quad q_\e\in \mP_{k-1}(\e), \e\in \cE_h(K),\\
	&\int_K\tauh:q_K\dx, \quad q_K\in [P_{k-2}(K)]^{2\times 2}_{\rm sym},\; K\in\cT_h.
\end{align}	

This leads to the construction of equilibrated moment tensor:
\begin{lemma} \label{Eqflux_C0ip}
	\cite[Lemma~5.1]{Bra_Pech_Sch_20_equi_flux_bih} There exists unique equilibrated moment tensor $\eqflux\in \Mh$ such that for each $K\in\cT_h$, 
	\begin{align}
		\eqfnn&=\avg{D^2 u_{\ip,nn}}_{\e}-\frac{\sigma}{h_\e}\jump{\partial_n u_{\ip}}_{\e} \quad\in \mP_{k-1}(\e), \e\in \cE_h(K),\\
		\int_K\eqflux:q_K\dx&=\int_KD^2 u_{\ip}:q_K\dx-\sum_{\e\subset\cE_h(K)}\int_\e\gamma_\e\jump{\partial_n u_{\ip}}_{\e}q_{K,nn}\dss\fl q_K\in  [P_{k-2}(K)]^{2\times 2}_{\rm sym},
	\end{align}
	where $\gamma_\e=1/2$ for interior edge $\e\in \cE_h(\Omega)$ and $\gamma_\e=1$ for a boundary edge $\e\in\cE_h(\partial\Omega)$. Moreover, the equilibrated moment tensor satisfies \cite[eq. (5.6)]{Bra_Pech_Sch_20_equi_flux_bih}
	\begin{equation}\label{eqn_equilib_c0ip}
		\langle\Div\Div \eqflux,v_{\ip}\rangle=(f,v_{\ip})\fl v_{\ip}\in V_h^k.
	\end{equation}
\end{lemma}	

By the above Lemma~\ref{Eqflux_C0ip} and following \cite{Bra_Pech_Sch_20_equi_flux_bih}, we have the efficiency result: 
\begin{lemma}\label{lmm:c0ip_eff}
	Let $u_{\ip}$ be the discrete solution of \eqref{c0ip_prim} and $\eqflux$ be of \eqref{eqn_equilib_c0ip}. Then the following efficiency result holds:
	\begin{equation}
		\|\eqflux-D^2u_{\ip}\|^2\lesssim\|u-u_{\ip}\|_{\ip}^2+\sum_{K\in\cT_h}h_K^4\|f-\bar{f}\|^2_{L^2(K)},
	\end{equation}	
	where $\bar{f}$ is any interpolation of $f$ into the space of piecewise polynomial functions of total degree less than equals to $k$.
\end{lemma}

\bigskip

\noindent\textbf{Computation of potential reconstruction.}
We describe potential reconstruction for $k=2$ by averaging \cite{BGS10,Ciarlet_78_Interpolation_P3HCT, Brenner}: let $\bar{V}_h\subset \hto$ be the Hsieh--Clough--Tocher associated with the triangulation $\cT_h$. For higher-degree approximations $k\geq 3$, we refer \cite{Georgoulis2011} for extension of this approach, see \eqref{dg_enrich} below. We define the enrichment operator $E_h: V_h^k\to\bar{V}_h$ as follows:  let $N$ be any (global) degree of freedom of $\bar{V}_h$, i.e., $N$ is either the evaluation of a shape function or its first-order derivatives at an interior vertex of $\cT_h$ or the evaluation of the normal derivative of a shape function at the midpoint of an interior edge. For $v_{\ip}\in V_h^k$  define
\begin{equation}\label{c0ip_potential}
	N(E_hv_{\ip})=\frac{1}{|\cT_N|}\sum_{K\in\cT_N} N(v_{\ip}|_K)
\end{equation}
where $\cT_N$ is the set of triangles in $\cT_h$ that share the degree of freedom $N$ and $|\cT_N|$ is the number of elements of $\cT_N$. 
The enrichment operator satisfies the estimate:
\begin{equation}\label{c0ip_enrich}
	\|E_h v_{\ip}-v_{\ip}\|_{\ip}\leq C \inf_{v\in\hto}\|v-v_{\ip}\|_{\ip},
\end{equation}
for some positive constant $C$. Finally, we set $s_h:=E_hu_{\ip}$ and $\ts_h:=E_h\tu_{\ip}$ to compute the estimator in \eqref{abs_goal_est}. Moreover, the efficiency $\|s_h-u_{\ip}\|_{\ip}\leq C\|u-u_{\ip}\|_{\ip}$ follows from \eqref{c0ip_enrich} with $v_{\ip}=u_{\ip}$ and a choice $v=u$.

\bigskip

\noindent\textbf{Computation of data oscillation.}  We follow the procedure of \cite[Lemma~6.1]{Bra_Pech_Sch_20_equi_flux_bih} to compute oscillation of data $f$ and $\tf$. Assume the data $f$ and $\tf$ belong to $L^2(\Omega)$. Let $\bar{f}$ denote the $L^2$ projection of $f$ onto the (discontinuous) space of piecewise polynomials of degree $k-3$ in $\cT_h$. Then the oscillation can be bounded by 
\begin{align}
	\|f-f_h\|_{-2}\leq c\left(\sum_{K\in\cT_h}h_K^4\|f-\bar{f}\|_{0,K}^2\right)^{1/2} \text{ and } \|\tf-\tf_h\|_{-2}\leq c\left(\sum_{K\in\cT_h}h_K^4\|\tf-\bar{\tf}\|_{0,K}^2\right)^{1/2},\label{c0ip_osc}
\end{align}
where the constant $c$ has upper bound $0.3682146$ due to an estimate of the interpolation by the Morley element, see \cite{Carst_Gal_eigs_bilap_14}. In the case of $k=2$, the projections can be set as $\bar{f}=0$ and $\bar{\tf}=0$.


We state convergence result, see  \cite{BS_C0IP05,BGS10,Gudi10}:
\begin{align}\label{IP_conv}
	\|u-u_{\ip}\|_{h}\leq \left(\inf_{v_{\ip}\in V_h^k}\|u-v_{\ip}\|_h+osc_{2}(f)\right)
\end{align}
where, the norm defined by  $\displaystyle\|v_{\ip}\|_{h}^2:=\|v_{\ip}\|_{\ip}^2+\sum_{\e\in\cE_h}\sum_{i,j=1,2}\left\|\avg{\frac{\partial^2 v_{\ip}}{\partial x_i\partial x_j}}_{\e}\right\|_{L^2(e)}^2$ and the data oscillation $\displaystyle osc_{2}(f):=\left(\sum_{K\in\cT_h}h_K^{4}\inf_{\bar{f}\in P_{k-2}(K)}\|f-\bar{f}\|^2_{L^2(K)}\right)^{1/2}$. This is used to obtain a convergence result for goal error as follows:

For the above $C^0$IPDG approximation, we observe that $\jump{u_\ip}_e=0=\jump{\partial_{\tau}u_\ip}_e$. This is used to simplify the goal residual estimator of \eqref{goal_main_est} as follows: 
\begin{align}
	\eta_h:&= (f-\Div\Div \eqflux,\tu_\ip)+ \sik(\eqflux-D^2u_\ip):\teqflux\dx+ \sum_{\e\in\cE_h}\int_e\jump{\partial_n u_\ip}_{\e}\teqfnn\dss.\label{simple_goal_est_IP}
\end{align}
Moreover, the remainder term has the estimate:
\begin{theorem}
	Let $u$ and $\tu$ respectively be the solution of \eqref{weak_prim} and \eqref{weak_dual}. Let $u_{\ip}\in V_h^k$ and $\tu_{\ip}\in V_h^m$ respectively be solution of \eqref{c0ip_prim} and \eqref{c0ip_dual}. Assume $\|u-u_{\ip}\|_{\ip}$ and  $\|\tu-\tu_{\ip}\|_{\ip}$, respectively, converge with orders $O(h^k)$ and $O(h^m)$. Then the remainder estimator $\cR_h$ of \eqref{goal_rem_est} has the convergence
	\begin{equation}\label{c0ip_pem_conv}
		|\cR_h(u,\tu,f;u_{\ip},\tu_{\ip})|\leq Ch^{k+m},
	\end{equation}
	where the positive constant $C$ (independent of the mesh parameter $h$) depends on load function $f$, and exact solutions $u$ and $\tu$.
\end{theorem}
\begin{proof}
	Recall the remainder estimator term $\cR_h(u,\tu,f;u_{\ip},\tu_{\ip})$ of \eqref{goal_rem_est}
	\begin{align}
		\cR_h(u,\tu,f;u_{\ip},\tu_{\ip})&=\langle f-f_h,\tu-\tsh\rangle + \sum_{K\in\cT_h}\int_K(\eqflux-D^2u_{\ip}):(D^2\tu-\teqflux)\dx\notag\\
		&\quad+\sum_{\e\in\cE_h}\int_e\jump{\partial_n u_{\ip}}_{\e}(D^2_{nn}\tu-\teqfnn)\dss.\label{c0ip_pem_est}
	\end{align}
	The first oscillation term in the above \eqref{c0ip_pem_est} is estimated by \eqref{osc_k} as
	\begin{equation}
		|\langle f-f_h,\tu-\tsh\rangle|\leq \|f-f_h\|_{-2}(\|\tf-\tf_h\|_{-2}+\|D^2 \ts_h-\teqflux\|).\label{est_osc_conv}
	\end{equation}
	The identity \eqref{eqnDual_prag_syn_osc} with $\tilde{v}=\tsh$ yields $\|D^2\tu-\teqflux\|\leq \|D^2\ts_{h}-\teqflux\|+\sqrt{2}\,osc_{\rm dual}(\tf,\tu)$ by Schwarz inequality. This leads to an estimate for  the second term in \eqref{c0ip_pem_est} as
	\begin{align}
		&\sum_{K\in\cT_h}\int_K(\eqflux-D^2u_{\ip}):(D^2\tu-\teqflux)\dx\leq \|\eqflux-D^2u_{\ip}\|\,\|D^2\tu-\teqflux\|\notag\\
		&\quad\leq \|\eqflux-D^2u_{\ip}\|\left(\|D^2\ts_{h}-\teqflux\|+\sqrt{2}\,osc_{\rm dual}(\tf,\tu)\right).
	\end{align}
	The last term of \eqref{c0ip_pem_est} is bounded by Cauchy--Schwarz inequality
	\begin{align}
		\Big{|}\sum_{\e\in\cE_h}\int_e\jump{\partial_n u_{\ip}}_{\e}(D^2_{nn}\tu-\teqfnn)\dss\Big{|}&\leq \sum_{e\in\cE_h} \|h_e^{-1/2}\jump{\partial_n u_{\ip}}_{\e}\|_{L^2(e)}\|h_e^{1/2}(D^2_{nn}\tu-\teqfnn)\|_{L^2(e)}\notag\\
		&\leq \|u-u_{\ip}\|_{\ip}\left(\sum_{\e\in\cE_h}\|h_e^{1/2}(D^2_{nn}\tu-\teqfnn)\|_{L^2(\e)}^2\right)^{1/2}.
	\end{align}
	Addition and subtraction of $u_{\ip}$ with trace inequality yield
	\begin{align}
		\sum_{\e\in\cE_h}\|h_e^{1/2}(D^2_{nn}\tu-\teqfnn)\|_{L^2(\e)}^2&\leq \sum_{\e\in\cE_h}\|h_e^{1/2}D^2_{nn}(\tu-\tu_{\ip})\|_{L^2(\e)}^2+\sum_{\e\in\cE_h}\|h_e^{1/2}(D^2_{nn}\tu_{\ip}-\teqfnn)\|_{L^2(\e)}^2\notag\\
		&\leq \|\tu-\tu_{\ip}\|_h+\|D^2\tu_{\ip}-\teqflux\|.\label{est_nrml}
	\end{align}
	The above estimates \eqref{est_osc_conv}-\eqref{est_nrml} and the efficiency result of Lemma~\ref{lmm:c0ip_eff} for primal and dual problems yield the required estimate \eqref{c0ip_pem_conv}.
\end{proof}

\begin{corollary}
	If the primal and dual solutions respectively $u$ and $\tu$ belong to $H^{2+\alpha}(\Omega)\cap\hto$, then the remainder estimator $\cR_h$ of \eqref{goal_rem_est} has the convergence
	\begin{equation}\label{c0ip_pem_low_peg}
		|\cR_h(u,\tu,f;u_{\ip},\tu_{\ip})|\leq Ch^{2\alpha},
	\end{equation}
	where the positive constant $C$ (independent of the mesh parameter $h$) depends on load function $f$, and exact solutions $u$ and $\tu$.
\end{corollary}

\subsection{Discontinuous Galerkin FEMs}

Let $V_h^k:=\mP_k(\cT_h)$. Define the bilinear form $a_{\dg}: V_h^k\times V_h^k\to \mR$  by \cite{Bra_Hop_Lin_DG_EqFlux_Bih_18}
\begin{align}
	a_{\dg}(u_{\dg},v_{\dg})&:=\sum_{K\in\cT_h}\int_K D^2 u_{\dg}:D^2 v_{\dg}\dx-\sum_{\e\in\cE_h}\int_\e \left(\jump{\nabla u_{\dg}}_{\e}{\cdot}\avg{D^2 v_{\dg}\bfn}_{\e}+\avg{D^2u_{\dg}\bfn}_{\e}{\cdot}\jump{\nabla v_{\dg}}_{\e}\right)\dss\notag\\
	&\qquad+\sum_{\e\in\cE_h}\int_\e \left(\jump{ u_{\dg}}_{\e}\avg{\Div(D^2 v_{\dg}){\cdot}\bfn}_{\e}+\avg{\Div(D^2u_{\dg}){\cdot}\bfn}_{\e}\jump{v_{\dg}}_{\e}\right)\dss\notag\\
	&\qquad\quad+\sum_{\e\in\cE_h}\frac{\sigma_1}{h_\e}\int_\e \jump{\partial_n u_{\dg}}_{\e}\jump{\partial_nv_{\dg}}_{\e}\dss+\sum_{\e\in\cE_h}\frac{\sigma_2}{h_\e^3}\int_\e \jump{ u_{\dg}}_{\e}\jump{v_{\dg}}_{\e}\dss
\end{align}
for positive penalty parameter $\sigma_1$ and $\sigma_2$, and linear forms
\begin{equation}
	l_{\dg}(v_{\dg}):=\sik fv_{\dg}\dx \text{ and } \tilde{l}_{\dg}(v_{\dg}):=\sik \tf v_{\dg}\dx.
\end{equation}
The DG method for \eqref{weak_prim} seeks $u_{\dg}\in V_h^k$ such that
\begin{equation}\label{dg_prim}
	a_{\dg}(u_{\dg},v_{\dg})=l_{\dg}(v_{\dg})\fl v_{\dg}\in V_h^k,
\end{equation}
and for the dual problem \eqref{weak_dual} seeks $\tu_{\dg}\in V_h^k$ such that
\begin{equation}\label{dg_dual}
	a_{\dg}(\tu_{\dg},v_{\dg})=\tilde{l}_{\dg}(v_{\dg})\fl v_{\dg}\in V_h^k.
\end{equation}
The discretization error will be measured by the mesh-dependent dG norm
\begin{equation*}
	\|v\|_{\rm dG}^2:=\sum_{K\in\cT_h}\|D^2 v\|_{0,K}^2+\sum_{\e\in\cE_h}\frac{\sigma_1}{h_\e}\|\jump{\partial_n v}_{\e}\|_{0,\e}^2+\sum_{\e\in\cE_h}\frac{\sigma_2}{h_\e^3}\|\jump{v}_{\e}\|_{0,\e}^2 \fl v\in V_h^k+\hto.
\end{equation*}
It is well known that for sufficiently large $\sigma_1=O((k+1)^2)$ and $\sigma_2=O((k+1)^6)$, there exists a positive constant $\beta$ such that the following coercivity
\begin{equation*}
	a_{\dg}(v_{\dg},v_{\dg})\geq \beta\|v_{\dg}\|_{\rm dG}^2 \fl v_{\dg}\in V_h^k,
\end{equation*}
and for all $ v_\dg,w_\dg\in V_{\dg}$, boundedness $|a_{\dg}(v_\dg,w_\dg)|\leq C\|v_\dg\|_{\rm dG}\|w_\dg\|_{\rm dG}$ result hold. Moreover, one can extend the definition of $a_{\dg}(\cdot,\cdot)$ to $V_h^k+\hto$ by a lifting operator, see \cite{Georgoulis2009}, and have the coercivity and boundedness of the extension. An abuse of notation, we also denote the extension of $a_{\dg}(\cdot,\cdot)$ to $V_h^k+\hto$ by itself.

\noindent\textbf{Construction of equilibrated moment tensor.}
We follow \cite{Bra_Hop_Lin_DG_EqFlux_Bih_18} for the construction of equilibrated moment tensor $\eqflux$ and $\teqflux$.
The equilibrated moment tensors are constructed in the discrete space $\Mh$ defined by
\begin{equation}
	\Mh:=\left\{\tauh\in L^2(\Omega)^{2\times 2}\; |\; \tauh|_K\in P_{\ell}(K)^{2\times 2} , K\in\cT_h\right\}\cap\Hdivsq,
\end{equation}
where $\ell:=\begin{cases}
	k\quad \text{if } k\geq 3,\\
	3\quad \text{if } k=2.
\end{cases}$

For $K\in\cT_h$, let $f_K$ be the $L^2$-projection of $f$ onto $P_{l-2}(K)$, and let $f_h\in\lt$ be such that $f_h|_K=f_K, K\in\cT_h$. Let $\BDM_{m}(K), m\in{\mathbb N}$ be denoted by the Brezzi--Douglas--Marini element of polynomial degree $m$, see \cite{Brez_For_91}.
The construction of equilibrated moment tensor is obtained by two steps: first construct an auxiliary vector field $\psibarheq\in\Hdiv,\; \psibarheq|_K\in\BDM_{\ell-1}(K), K\in\cT_h$ satisfying
\begin{equation}
	\nabla{\cdot}\psibarheq=f_h\quad\text{ in } L^2(\Omega),
\end{equation}
and then an equilibrated moment tensor $\eqflux\in\Mh$ satisfying
\begin{equation}
	\nabla{\cdot}\eqflux=\psibarheq\quad\text{ in } L^2(\Omega).
\end{equation}
Define some auxiliary numerical flux functions on the edges $\e\in\cE_h$ by
\begin{align*}
	\underline{\widehat{\boldsymbol{u}}}^{(1)}&:=
	\begin{cases}
		\left\{\hskip -4pt\left\{\nabla u_h\right\}\hskip -4pt\right\}_\e,\; \e\in\cE_h(\Omega)\\
		\qquad\quad\; 0,\; \e\in\cE_h(\partial\Omega),
	\end{cases}\\
	\underline{\widehat{u}}^{(2)}&:=
	\begin{cases}
		\left\{\hskip -4pt\left\{u_h\right\}\hskip -4pt\right\}_\e,\; \e\in \cE_h(\Omega)\\
		\qquad\;\; 0,\; \e\in\cE_h(\partial\Omega),
	\end{cases}\\
	\underline{\underline{\widehat{\boldsymbol{p}}}}&:=\avg{D^2 u_h}_\e-\frac{\sigma_1}{h_\e}\bfn\jump{\nabla u_h}_\e^T,\\
	\psibarhat&:=\avg{\nabla{\cdot} D^2 u_h}_\e+\frac{\sigma_2}{h_\e^3} \jump{u_h}_\e\bfn.
\end{align*}

The auxiliary vector field $\psibarheq$ is constructed locally on each element $K\in\cT_h$ such that $\psibarheq\in \BDM_{\ell-1}(K)$ satisfies the following interpolation conditions, see \cite[Eq. 6.5]{Bra_Hop_Lin_DG_EqFlux_Bih_18}
\begin{align*}
	\int_{\e}\bfn \cdot \psibarheq\, q\dss&=\int_{\e}\bfn\cdot\psibarhat\, q\dss, \quad q\in P_{\ell-1}(\e), \; \e\in\cE_h(\partial K),\\
	\int_{K}\psibarheq \cdot \nabla q\dx&=\int_{\partial K} \bn_{\partial K}\cdot\psibarhat q\dss-\int_K fq\dx, \quad q\in P_{\ell-2}(K),\\
	\int_K\psibarheq \cdot\curl(b_K q)\dx&=\int_K(\nabla{\cdot}D^2u_h) \cdot \curl(b_Kq)\dx, \quad
	q\in P_{\ell-3}(K),
\end{align*}
where $b_K=\lambda_1^K\lambda_2^K\lambda_3^K$ is the bubble function on element $K$ for barycentric coordinates $\lambda_i^K, i=1,2,3$ of $K$. Finally, the equilibrated moment tensor $\eqflux=(\sigma_{ij}^{h,{\rm eq}})_{i,j=1}^2$, with $\underline{\boldsymbol{\sigma}}_{h,{\rm eq}}^{(i)}=(\sigma_{i1}^{h,{\rm eq}},\sigma_{i2}^{h,{\rm eq}})^T,\; 1\leq i\leq 2$, in each element $K$ is constructed by fixing the degrees of freedom \cite[Eq. 6.8]{Bra_Hop_Lin_DG_EqFlux_Bih_18}:
\begin{align*}
	\int_{\e}\eqflux\bfn \cdot \underline{\boldsymbol{q}}\dss&=\int_{\e} \underline{\underline{\widehat{\boldsymbol{p}}}}\bfn \cdot \underline{\boldsymbol{q}}\dss, \quad\underline{\boldsymbol{q}}\in [P_{\ell}(\e)]^2, \; \e\in\cE_h(\partial K),\\
	\int_{K}\eqflux : \nabla \underline{\boldsymbol{q}}\dx&=-\int_{K}\psibarheq \cdot \underline{\boldsymbol{q}}\dx+\int_{\partial K}\underline{\underline{\widehat{\boldsymbol{p}}}}\bn_{\partial K} \cdot \underline{\boldsymbol{q}}\dx, \quad\underline{\boldsymbol{q}}\in [P_{\ell-1}(K)]^2,\\
	\int_K\underline{\boldsymbol{\sigma}}_{h,{\rm eq}}^{(i)} \cdot \curl(b_Kq)\dx&=\int_K \underline{\boldsymbol{z}}^{(i)} \cdot \curl(b_Kq)\dx, \quad q\in P_{\ell-2}(K), 1\leq i\leq 2,
\end{align*}
where $\underline{\boldsymbol{z}}^{(i)}=(\frac{\partial^2 u_h}{\partial x_i\partial x_1}, \frac{\partial^2 u_h}{\partial x_i \partial x_2}), i=1,2$. The above constructions lead to the equilibrium: $\Div\Div \eqflux=f_h$. A similar construction for dual problem with data $\tf$ and approximation $\tu_h$ leads to the equilibrated moment tensor $\teqflux$.

\noindent\textbf{Computation of potential reconstruction.} Let $S_h^r$ be a $C^1$-conforming finite-element space consisting of the macro-elements of order $r\geq 4$, see \cite[Definition 3.1]{Georgoulis2011}. We follow the construction of recovery operator of \cite{Georgoulis2011}. For each nodal point $\nu$ of the $C^1$-conforming finite-element space $S_h^{k+2}$, define $\omega_\nu$ to be the set of $ K\in\cT_h$ that share the nodal point $\nu$, i.e., $\omega_\nu=\{ K\in\cT_h: \nu\in K\}$. Define the operator $E_h: \mP_k(\cT_h)\to S_h^{k+2}\cap \hto$ by the averaging:
\begin{equation}\label{dg_enrichment}
	N_{\nu}(E_hv_{\dg})=
	\begin{cases}
		\frac{1}{|\omega_{\nu}|}\sum_{ K \in\omega_{\nu}} N_{\nu}(v_{\dg}|_{ K}) \quad\text{if } \nu\notin\partial\Omega,\\
		0\;\;\quad\qquad\qquad\qquad\qquad\text{if } \nu\in\partial\Omega,
	\end{cases}
\end{equation}
where $N_{\nu}$ is any nodal variable at $\nu$ and $\nu$ is any nodal point of $S_h^{k+2}$.
This operator satisfies the estimate \cite[Lemma 3.1]{Georgoulis2011} and \cite[Lemma 3.5]{CC_GM_NN_19}:
\begin{equation}\label{dg_enrich}
	\|E_h v_{\dg}-v_{\dg}\|_{\dg}\leq C \inf_{v\in\hto}\|v-v_{\dg}\|_{\dg},
\end{equation}
for some positive constant $C$.


\section{Numerical Experiments}\label{sec:numer}
In this section, numerical results of goal-oriented a posteriori estimations are presented for the $C^0$IPDG method of Section~\ref{sec:c0ip} with $k=2$. The approximate goal functional is defined by
\begin{equation}
	Q_h:=Q(u_h)+\left(\eqflux-D^2 s_h,\tbfsigh^{\rm m}\right).
\end{equation}
The primal and dual estimators are defined respectively by $\eta_h:=\|D^2s_h-\eqflux\|$ and $\tilde{\eta}_h:=\|D^2\ts_h-\teqflux\|$. This gives the following  error estimate from \eqref{abs_goal_est}
\begin{align}\label{num_goal_err_est}
	e_{h,{\rm goal}}:=|Q(u)-Q_h|\leq \frac{\eta_h\tilde{\eta}_h}{2}+|Q(s_h-u_h)|=:\eta_{h,{\rm goal}}^{\rm abs},
\end{align}
where the higher order data oscillation terms are not considered in the computations.
The estimators are further localized for a mesh adaptation as 
\begin{align}
	\eta_h^2&=\sum_{K\in\cT_h}\eta_{h,K}^2 \text{ where } \eta_{h,K}:=\|D^2s_h-\eqflux\|_{L^2(K)},\\
	\tilde{\eta}_h^2&=\sum_{K\in\cT_h}\tilde{\eta}_{h,K}^2 \text{ where } \tilde{\eta}_{h,K}:=\|D^2\ts_h-\teqflux\|_{L^2(K)}\text{ and }\\
	\eta_{h,{\rm NC}}^2&=\sum_{K\in\cT_h}\eta_{h,K,{\rm NC}}^2 \text{ where } \eta_{h,K,{\rm NC}}:=|Q((s_h-u_h)\chi_K)|,
\end{align}
and $\chi_K$ is the characteristic function defined on $K\in\cT_h.$
We apply Algorithm~\ref{alg:adap} which follows standard adaptive procedure  SOLVE, ESTIMATE, MARK and REFINE for the numerical examples below. For the experiments below, the penalty parameter $\sigma$ for the $C^0$IPDG method is set to $20$.

\begin{algorithm}
	\caption{Goal-oriented adaptive method} \label{alg:adap}
	Input: Initial mesh $\cT_0$, $J \geq 1$ the maximum number of mesh refinement levels, and two real parameters $\theta \in (0,1)$.
	
	\medspace
	
	Set $j=0$.
	\\ While ($j \leq J$) do
	\begin{itemize}
		\item {\bf SOLVE/COMPUTE}:
		\begin{enumerate}
			\item \label{item_1}  Solve the {\it primal and dual matrix systems} $\mA {\rm U}_j=\alg{F}_j$ and $\mA\tilde{{\rm U}}_j=\tilde{{\rm F}}_j$ related to the discrete problems.
			\item  Compute the potential reconstructions for primal $s_j$ and for dual $\ts_j$ of Definition~\ref{Defn_Potn_pecons}.  Compute the moment tensor for primal problem $\underline{\underline{\boldsymbol \sigma}}_j^{\rm eq}$ and for dual problem $\underline{\underline{\tilde{\boldsymbol \sigma}}}_j^{\rm eq}$ from Definition~\ref{Defn_moment_tensor}.		
		\end{enumerate}
		\item {\bf ESTIMATE}. Compute the primal estimator $\eta_{j}$,  the dual estimator $\tilde{\eta}_{j}$ and the nonconforming estimator $\eta_{j,{\rm NC}}$ proposed for goal-oriented error estimation.
		\item {\bf MARK}. Mark sets for each of the primal and dual problems:
		\begin{enumerate}
			\item The D\"{o}rfler marking chooses a minimal subset  $\cM_{j}^{\rm p}\subset \cT_{j}$ such that
			\begin{align*}
				\theta\, \sum_{K\in\mathcal{T}_{j}}\eta^2_{j}(K)\leq \sum_{K\in\mathcal{M}_{j}^{\rm p}}\eta^2_{j}(K).
			\end{align*}
			\item The D\"{o}rfler marking chooses a minimal subset  $\cM_{j}^{\rm d}\subset \cT_{j}$ such that
			\begin{align*}
				\theta\, \sum_{K\in\mathcal{T}_{j}}\tilde{\eta}^2_{j}(K)\leq \sum_{K\in\mathcal{M}_{j}^{\rm d}}\tilde{\eta}^2_{j}(K).
			\end{align*}
			\item The D\"{o}rfler marking chooses a minimal subset  $\cM_{j}^{\rm NC}\subset \cT_{j}$ such that
			\begin{align*}
				\theta\, \sum_{K\in\mathcal{T}_{j}}\eta^2_{j,{\rm NC}}(K)\leq \sum_{K\in\mathcal{M}_{j}^{\rm NC}}\eta^2_{j,{\rm NC}}(K).
			\end{align*}
			\item Set $\cM_{j}:=\cM_{j}^{\rm p}\cup\cM_{j}^{\rm d}\cup\cM_{j}^{\rm NC}$ the union of marked sets found for primal, dual and nonconforming marking procedures above.
		\end{enumerate}
		\item {\bf REFINE}.  Compute the closure of $\cM_{j}$ and generate a new triangulation $\cT_{j+1}$ using newest vertex bisection method (\cite{Stevenson08}). 
		
		Set $j:=j+1$.
	\end{itemize}
	End While
\end{algorithm}

In the following numerical tests, we also compute the estimator found in \eqref{simple_goal_est_IP} in the context of $C^0$IPDG method. We rename it as goal residual estimator
\begin{align}
	\eta_{h,{\rm goal}}^{\rm res}:&= \sik(\eqflux-D^2u_{\ip}):\teqflux\dx + \sum_{\e\in\cE_h}\int_e\jump{\partial_n u_{\ip}}_{\e}\teqfnn\dss.\label{goal_est_res}
\end{align}
The potential reconstructions $s_h$ and $\ts_h$ for primal and dual solutions are computed from the definition in \eqref{c0ip_potential}. The symmetric piecewise linear equilibrated moment tensors $\eqflux$ for primal and $\teqflux$ for dual problems are constructed from Lemma~\ref{Eqflux_C0ip}. The effectivity indices are computed by the ratio $\eta_{h,{\rm goal}}^{\rm abs}/e_{h,{\rm goal}}$ for abstract goal estimator $\eta_{h,{\rm goal}}^{\rm abs}$ and by $\eta_{h,{\rm goal}}^{\rm res}/e_{h,{\rm goal}}$ for residual type goal estimator $\eta_{h,{\rm goal}}^{\rm res}$.


\subsection{Regular solution and uniform refinements}\label{example_1}
\begin{figure}
	\begin{center}
		\includegraphics[width=0.4\textwidth]{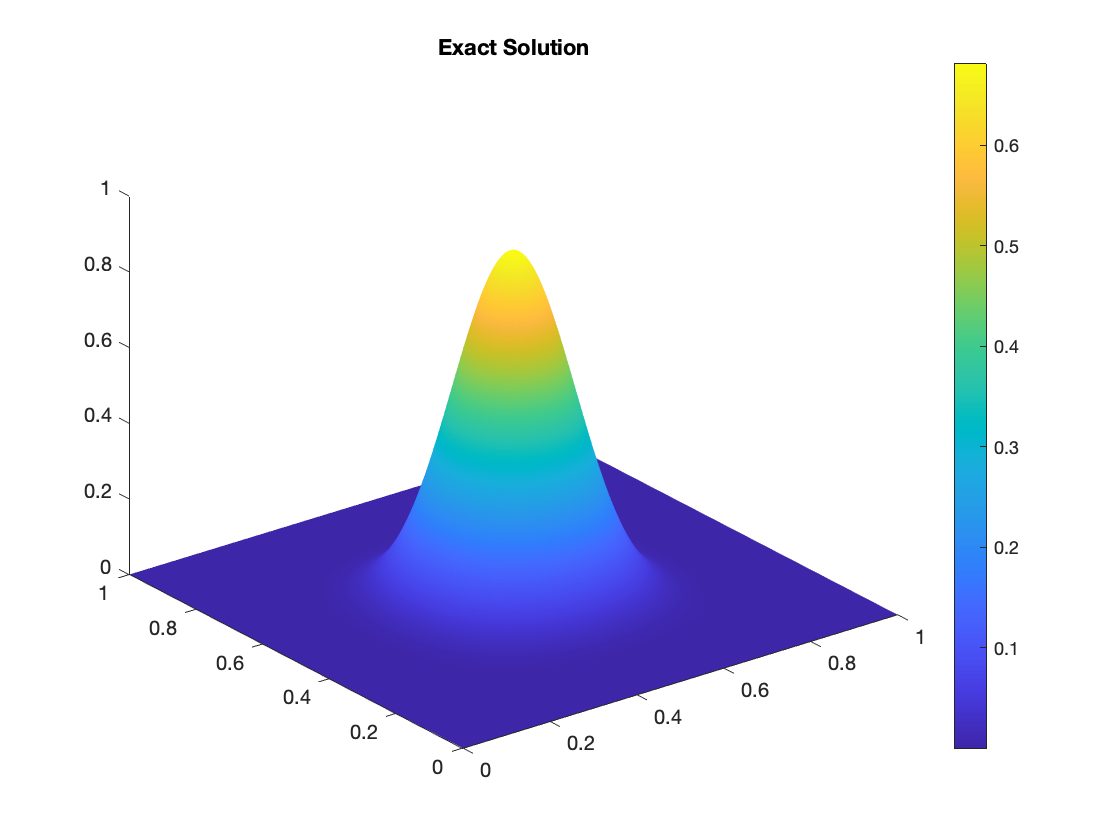}
		\begin{tikzpicture}[scale=3.5]
			\newcommand{\outline}[0]{%
				plot[smooth cycle]%
				coordinates{%
					(0.5,1),(1,0.5),(1,0.75),0.75,1)%
			}}
			\draw (0,0)--(1,0)--(1,1)--(0,1)--(0,0);
			\fill[gray] (0,1)--(0,0.75)--(0.75,0)--(1,0)--(1,0.25)--(0.25,1)--cycle;
			\draw (0.5,0.5) node {\textbullet};
			\node at (0.08,-.07) {$(0,0)$};
			\node at (0.9,-.07) {$(1,0)$};
			\node at (0.08,1.07) {$(0,1)$};
			\node at (0.4,0.6) {$\omega$};
		\end{tikzpicture}
		\caption{Exact solution ({\em left}) and the zone of interest ({\em right}). Example~\ref{example_1}, goal functional~\eqref{defn_goal_strip}.}
		\label{fig:exact_exp_zone}
	\end{center}
\end{figure}
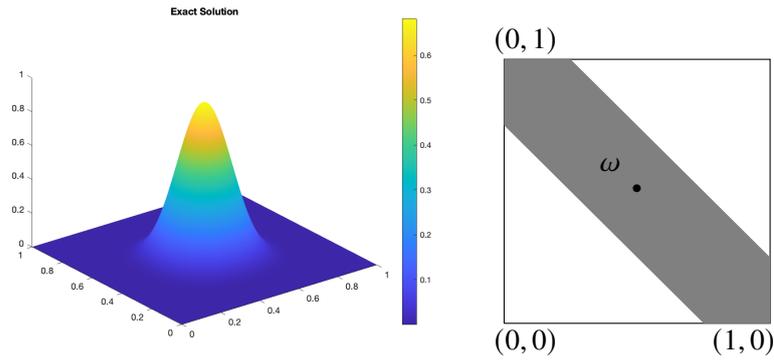

In this test, we consider an exact solution defined on a plate $\Omega:=(0,1)\times(0,1)$
\begin{equation}
	u(x,y)=10^{12}x^{10}(1-x)^{10}y^{10}(1-y)^{10}
\end{equation}
with load function $f$ defined by $f:=\Delta^2 u$ in $\Omega$. 
In other words, the goal functional is the mean value of the deflection in the strip $\omega$, where the right-hand side function $f$, the solution $u$ and gradient of $u$ exhibit large changes.
The exact solution has been illustrated in the left part of Figure~\ref{fig:exact_exp_zone} and the zone of interest $\omega$ is highlighted by grey color in the right part of Figure~\ref{fig:exact_exp_zone}. The solution has peak at $(\half,\half)$ which is highlighted by a bullet. The goal functional is defined by 
\begin{equation}\label{defn_goal_strip}
	Q(u)=\frac{1}{|\omega|}\int_\omega u\dx=\big(\tf,u\big)_{\Omega}, \text{ with } \tf = \frac{\chi_{\omega}}{|\omega|},
\end{equation}
where the strip $\omega:=\{(x,y)\in\Omega: 0.75\leq x+y\leq 1.25 \}$ illustrated in the right side of Figure~\ref{fig:exact_exp_zone} and $\chi_{\omega}$ is the characteristic function defined on $\omega$. The numerical integration value of the exact goal functional reads $Q(u)\approx 0.06044290015$. 

\begin{figure}
	\begin{center}
		\subfloat[]{\includegraphics[width=0.3\textwidth]{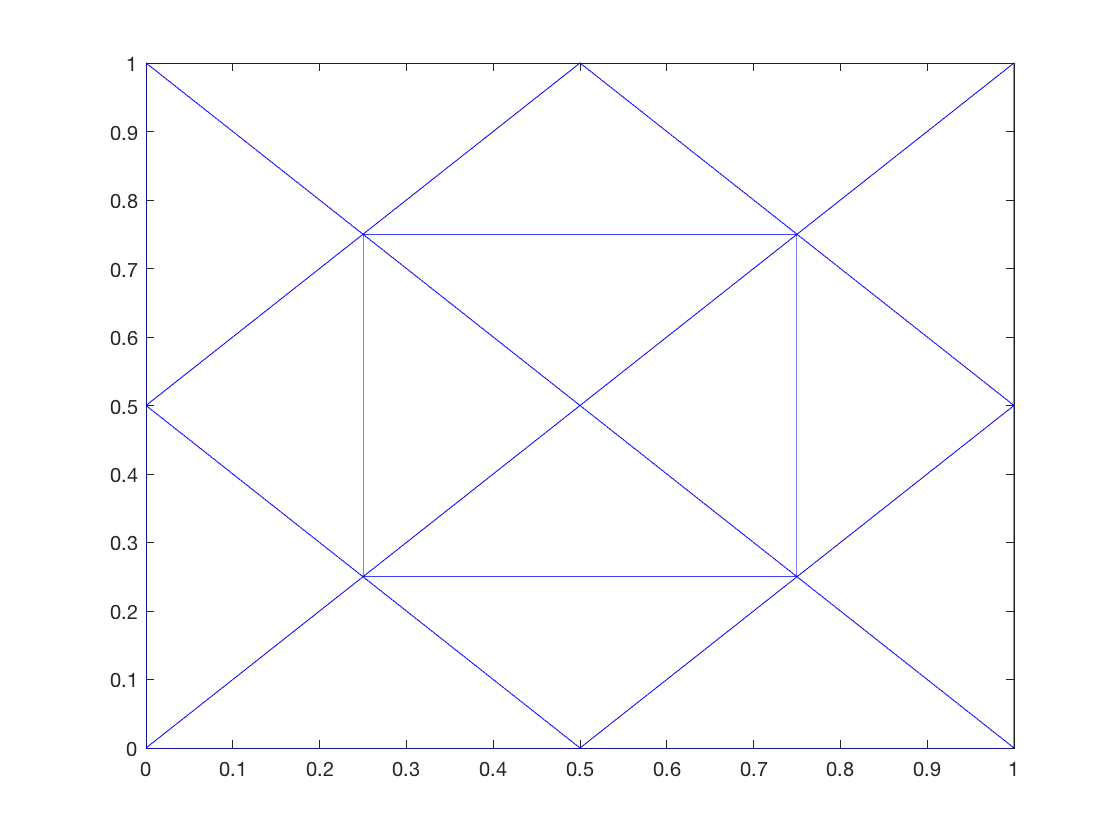}}
		\subfloat[]{\includegraphics[width=0.3\textwidth]{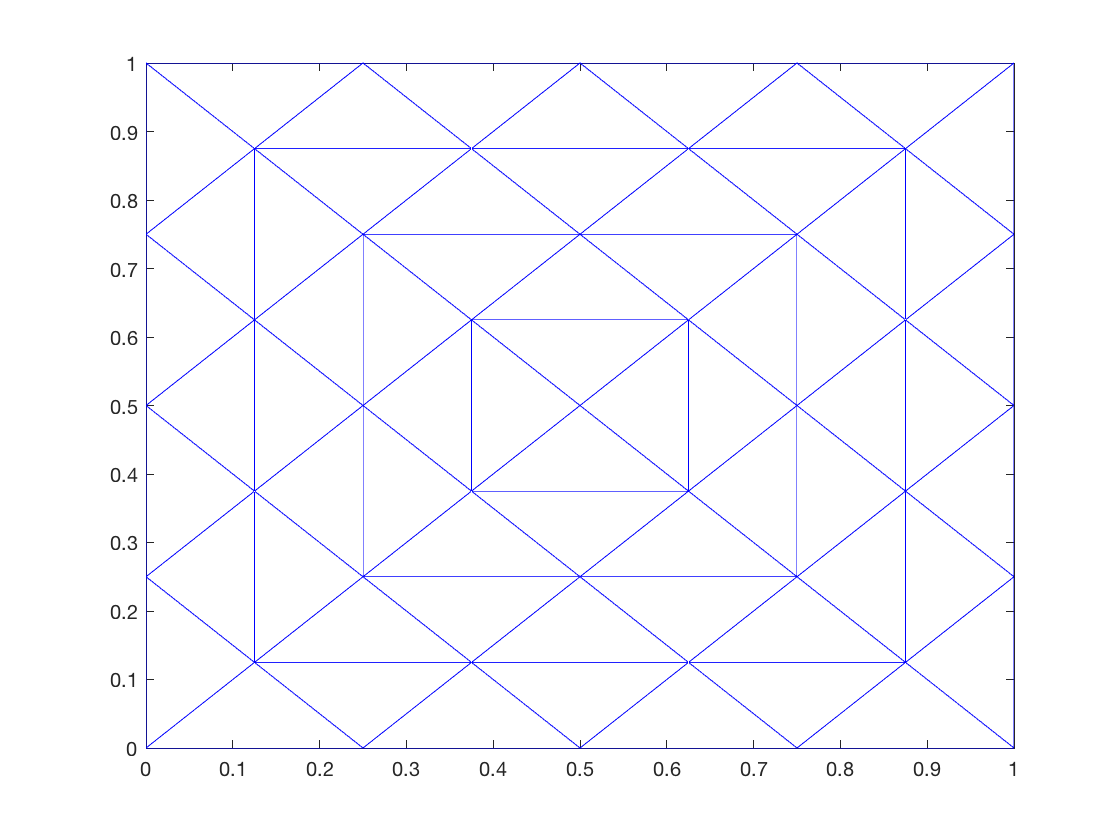}}
	\end{center}
	\caption{(a) Initial triangulation $\cT_0$ and (b) first uniform refinement $\cT_1$ of Example~\ref{example_1}.}
	\label{fig:T0}
\end{figure}

\begin{figure}
	\centering
	\includegraphics[height=0.24\textwidth]{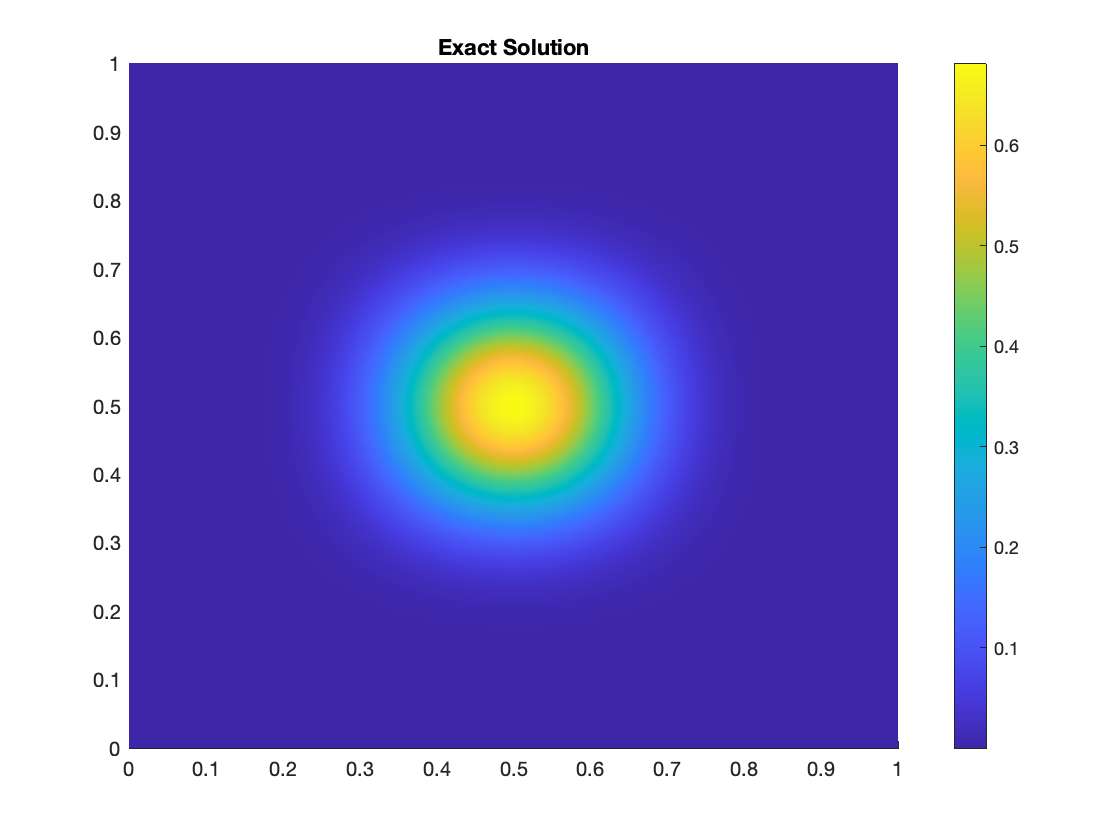}
	\includegraphics[height=0.24\textwidth]{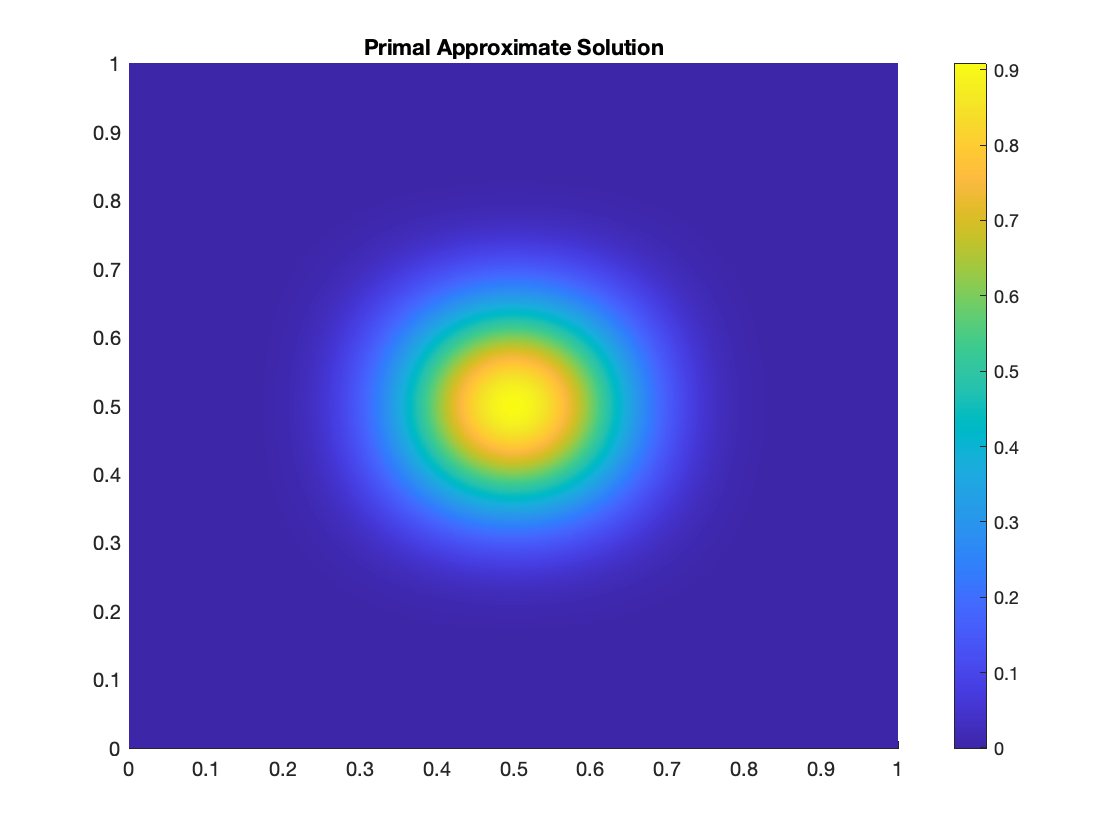}
	\includegraphics[height=0.24\textwidth]{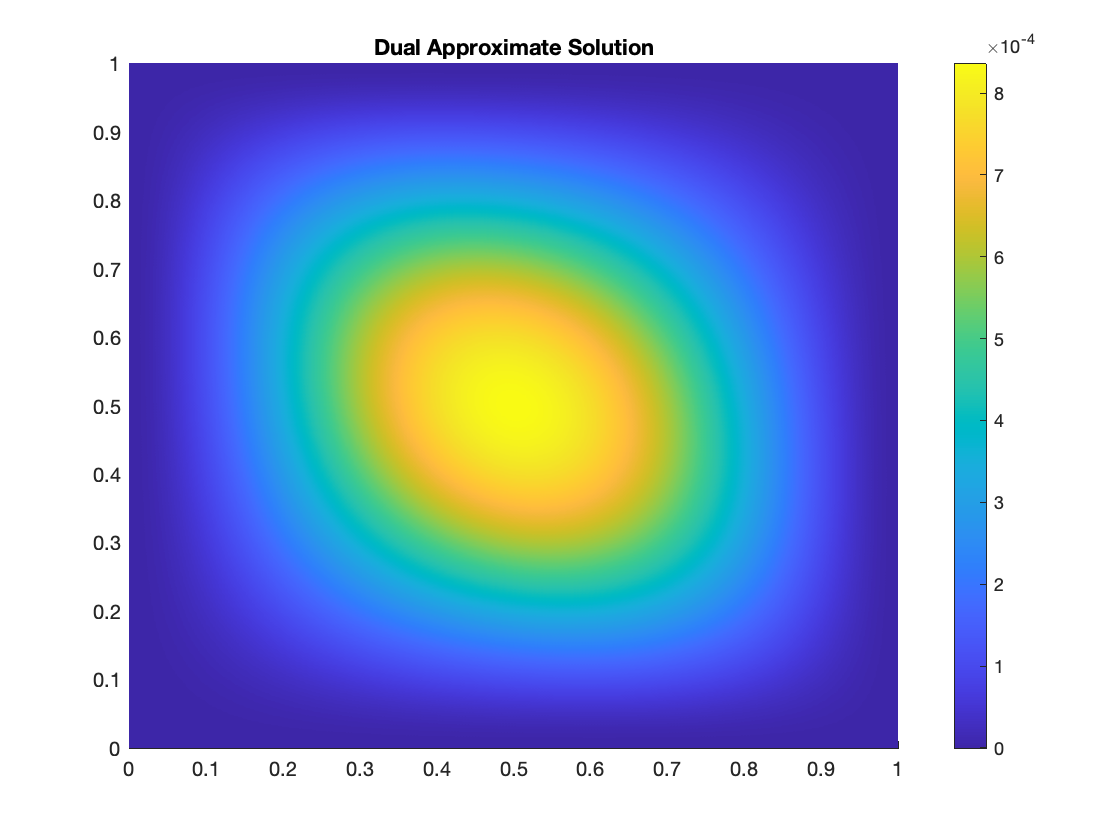}
	\caption{Exact primal solution ({\em left}), approximate primal solution ({\em middle}), and approximate dual solution ({\em right}). Example~\ref{example_1}, goal functional~\eqref{defn_goal_strip}.}
	\label{fig:primal_dual_apx}
\end{figure}

\begin{figure}
	\begin{center}
		\includegraphics[width=0.6\textwidth]{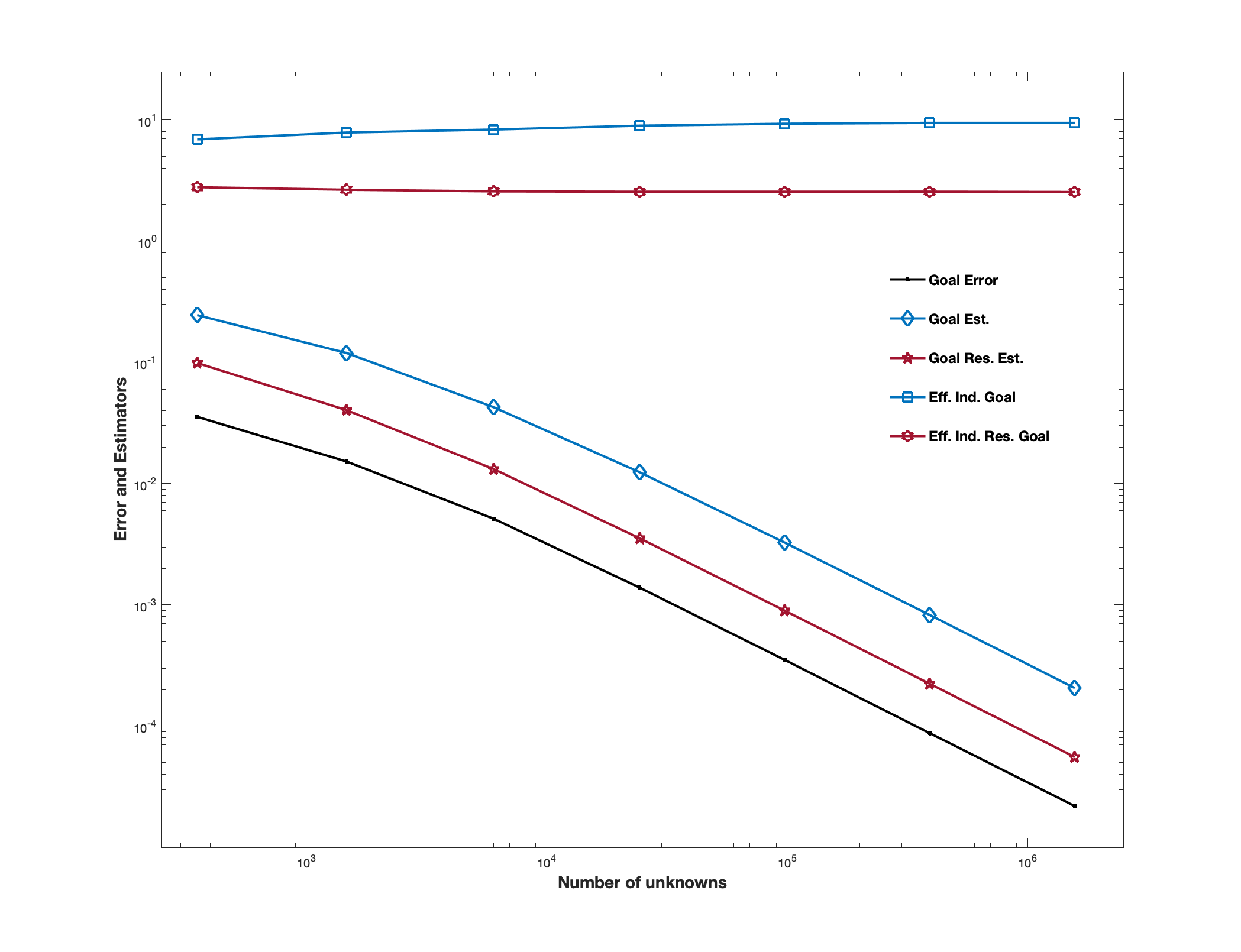}
	\end{center}
	\caption{Convergence histories of goal error $e_{h,{\rm goal}}$, abstract goal estimator $\eta_{h,{\rm goal}}^{\rm abs}$ and residual estimator $\eta_{h,{\rm goal}}^{\rm res}$ with effectivity indices.  Example~\ref{example_1}, goal functional~\eqref{defn_goal_strip}.}
	\label{fig:conv_hist_Sq}
\end{figure}

Numerical experiments are performed on the sequence of uniform  triangulations $\cT_0,\cT_1,\ldots,\cT_5$ with initial triangulation shown in Figure~\ref{fig:T0}(a). In the uniform refinement process each triangle is subdivided into four similar triangles, see Figure~\ref{fig:T0}. In Figure~\ref{fig:primal_dual_apx}, the exact solution $u$ in the left, approximate primal solution $u_{\ip}$ in the middle and approximate dual solution $\tu_{\ip}$ in the right are projected on the domain $\Omega$. The approximation for goal function is found to be $Q_h=0.06046477792$ on the mesh $\cT_5$. The convergence histories for goal error, goal estimator of \eqref{num_goal_err_est} and \eqref{goal_est_res} with respect to number of unknowns are plotted in Figure~\ref{fig:conv_hist_Sq}. We observe the quadratic convergence rates for goal error and goal estimators with effectivity index close to $9.4$ for abstract goal estimator \eqref{num_goal_err_est} and $2.5$ for \eqref{goal_est_res}.

\subsection{Singular solution and adaptive mesh refinement}\label{example_2}
\begin{figure}
	\begin{center}
		\includegraphics[width=0.5\textwidth]{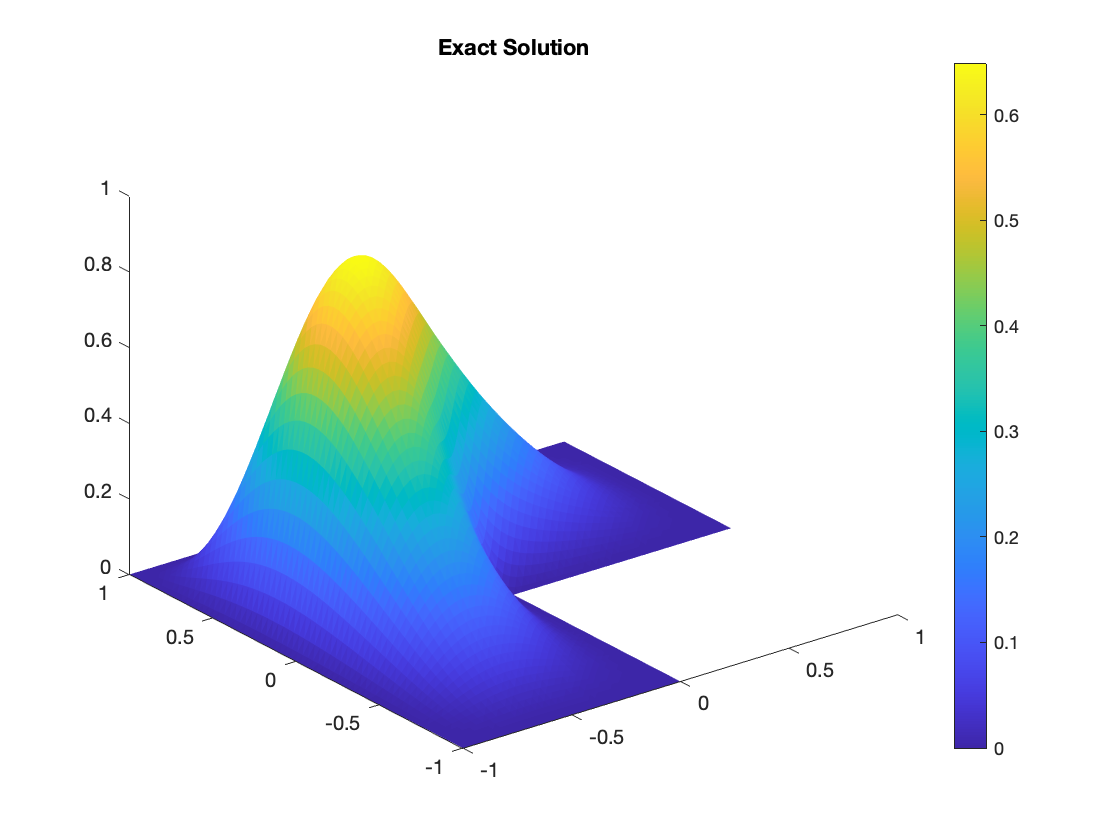}
		\hspace{0.25cm}
		\begin{tikzpicture}[scale=1.9]
			\draw (-1,-1)--(0,-1)--(0,0)--(1,0)--(1,1)--(-1,1)--(-1,-1);
			\draw (0.25,0) arc (0:270:0.25);
			\node at (0,0.1) {$\omega$};
			\node at (-0.5,0.5) {$\Omega$};
		\end{tikzpicture}
		\caption{Exact solution ({\em left}) and the zone of interest ({\em right}). Example~\ref{example_2}, goal functional~\eqref{defn_goal_sing_zone}.}
		\label{fig:Lshaped_exact_zone}
	\end{center}
	
\end{figure}
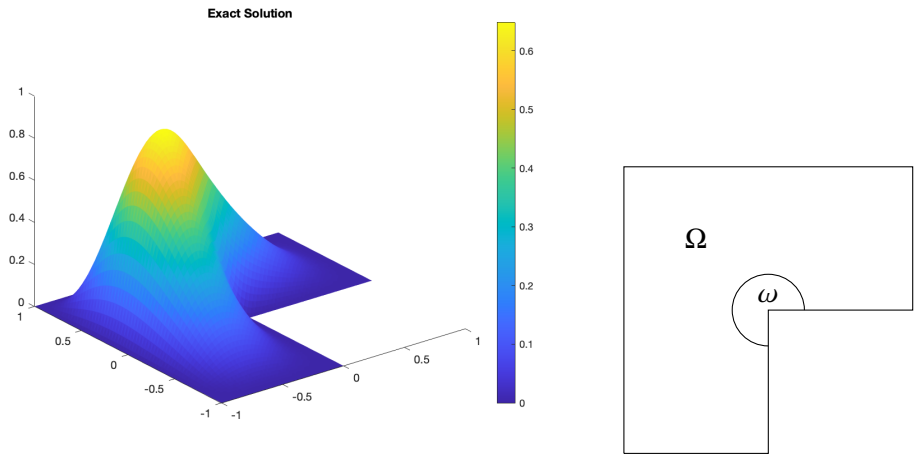

In this test, we consider the L-shaped domain  $\Omega=(-1,1)^2 \setminus\big{(}[0,1)\times(-1,0]\big{)}$. Set the
singular functions \cite{Grisvard}
$\displaystyle
u(r,\theta):=(1-r^2 \cos^2\theta)^2 (1-r^2 \sin^2\theta)^2 r^{1+\alpha}g_{\alpha,\omega}(\theta)$ 	with
$	g_{\alpha,\omega}(\theta):=$
\begin{align*}
	&\left(\frac{1}{\alpha-1}\sin\big{(}(\alpha-1)\omega\big{)}-\frac{1}{\alpha+1}\sin\big{(}(\alpha+1)\omega\big{)}\right)\times\Big{(}\cos\big{(}(\alpha-1)\theta\big{)}-\cos\big{(}(\alpha+1)\theta\big{)}\Big{)}\\
	&-\left(\frac{1}{\alpha-1}\sin\big{(}(\alpha-1)\theta\big{)}-\frac{1}{\alpha+1}\sin\big{(}(\alpha+1)\theta\big{)}\right)\times\Big{(}\cos\big{(}(\alpha-1)\omega\big{)}-\cos\big{(}(\alpha+1)\omega\big{)}\Big{)},
\end{align*}
where the angle $\omega:=\frac{3\pi}{2}$ and the parameter $\alpha= 0.5444837367$ is a non-characteristic root of $\sin^2(\alpha\omega) = \alpha^2\sin^2(\omega)$. It can be observed that the solution has the regularity $H^{2+\alpha}(\Omega)\cap \hto$, see \cite{Grisvard}. Since the problem has singularity at the origin $(0,0)$, we consider the goal functional
\begin{equation}\label{defn_goal_sing_zone}
	Q(u)=\frac{1}{|\omega|}\int_\omega u\dx=\big(\tf,u\big)_{\Omega}, \text{ with } \tf = \frac{\chi_{\omega}}{|\omega|},
\end{equation}
where $\omega:=\{(x,y)\in\Omega: (x-0)^2+(y-0)^2\leq 0.25^2 \}$ and $\chi_{\omega}$ is the characteristic function defined on $\omega$.
The exact solution ({\it left}), the domain $\Omega$ and the zone of interest ({\it right}) are illustrated in Figure~\ref{fig:Lshaped_exact_zone}. The numerical integration value of the exact goal functional reads $Q(u)\approx 0.018334438
$. 

\begin{figure}
	\begin{center}
		\subfloat[]{\includegraphics[width=0.25\textwidth]{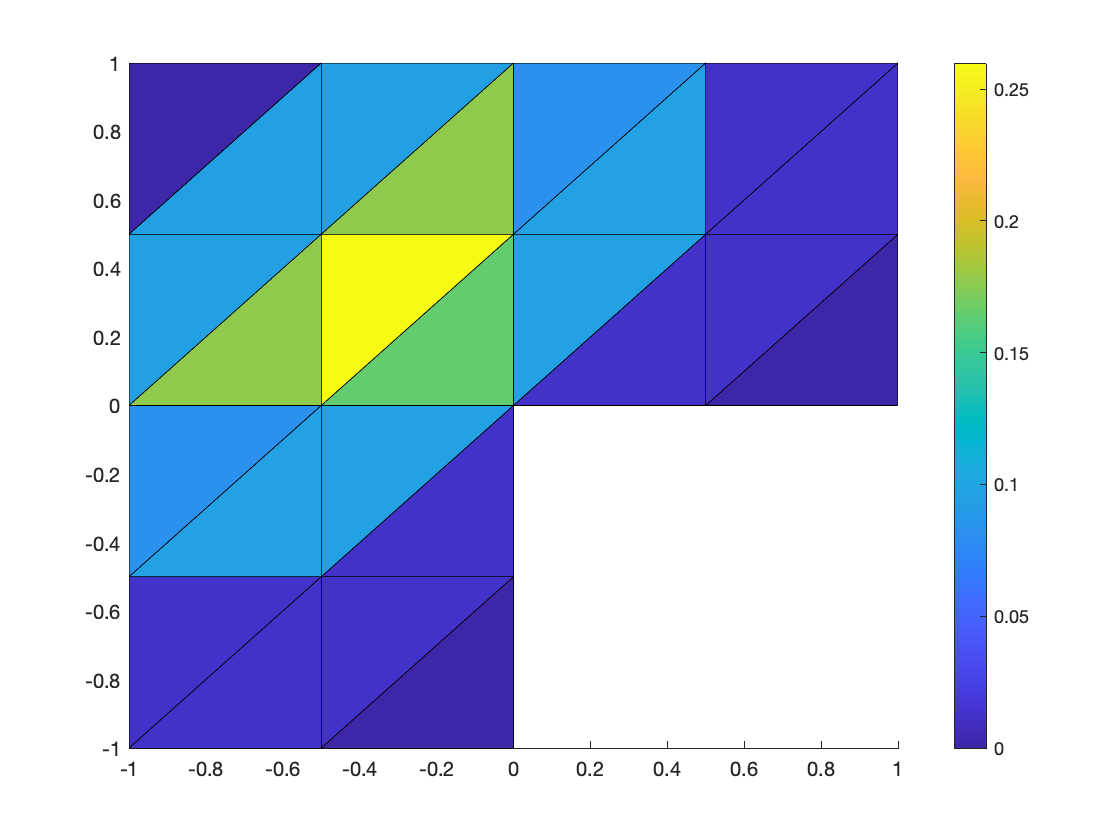}}
		\subfloat[]{\includegraphics[width=0.25\textwidth]{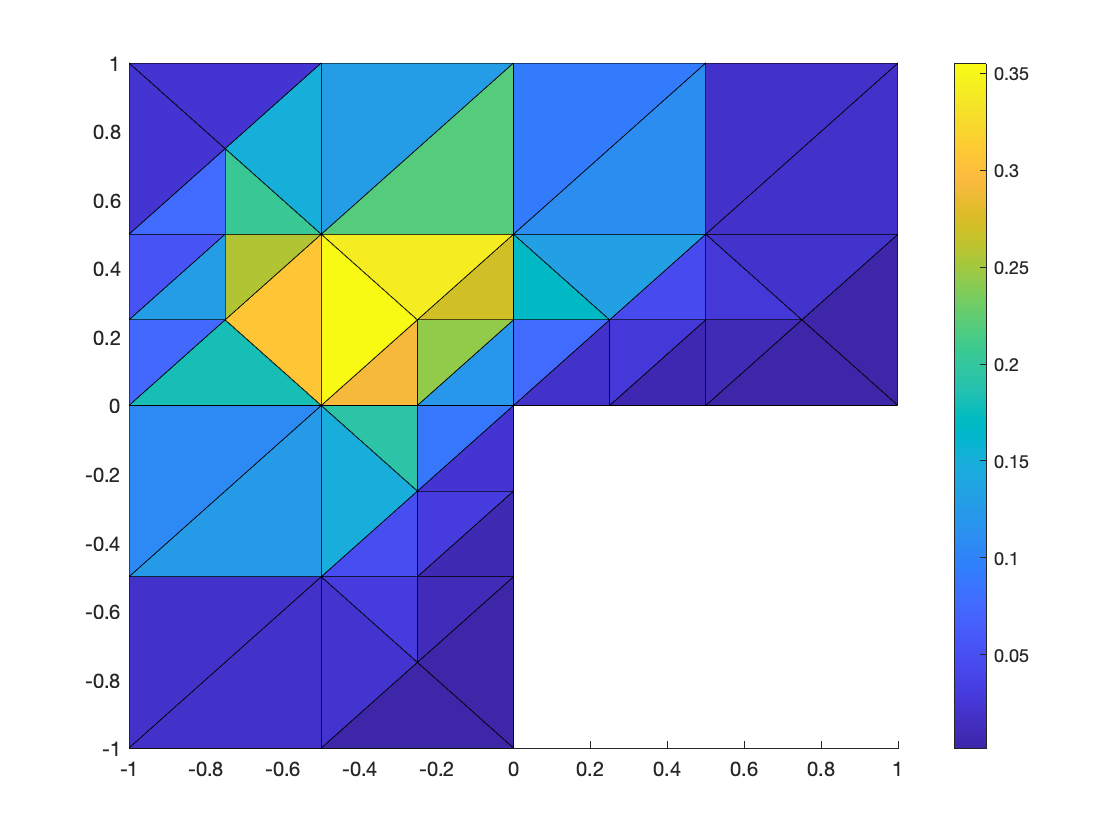}}
		\subfloat[]{\includegraphics[width=0.25\textwidth]{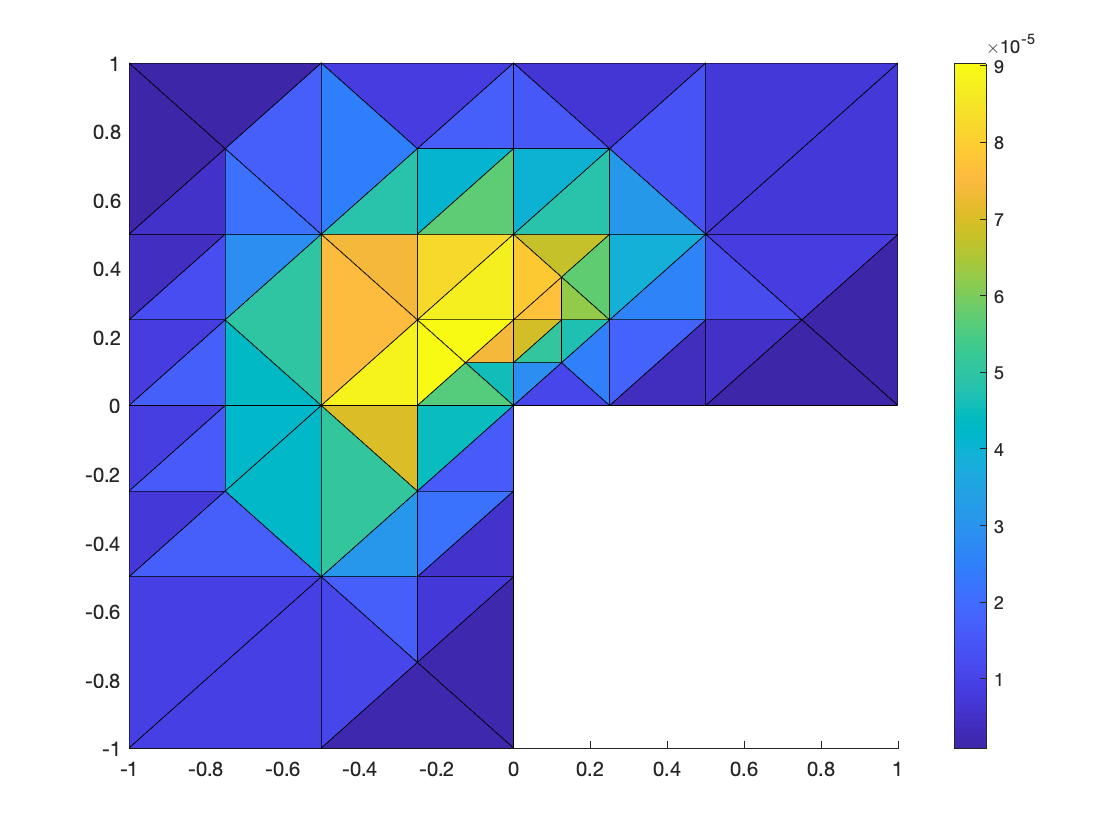}}
		\subfloat[]{\includegraphics[width=0.25\textwidth]{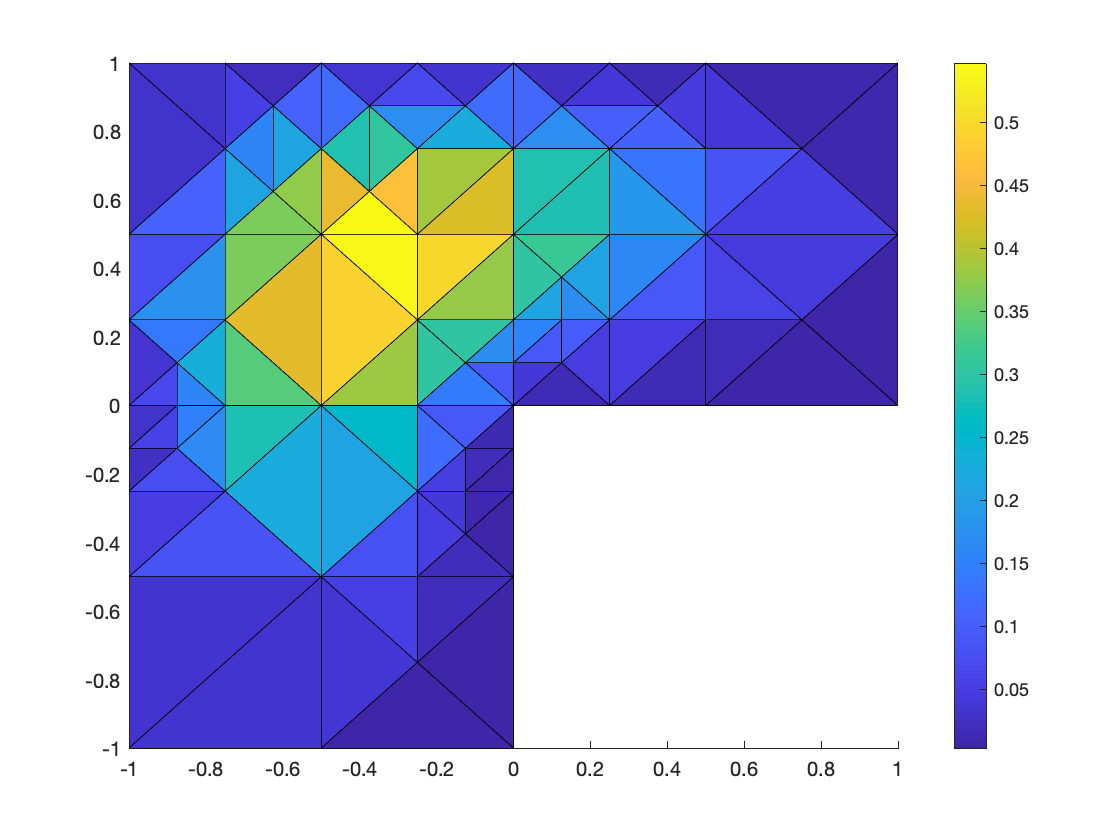}}\\
		\subfloat[]{\includegraphics[width=0.25\textwidth]{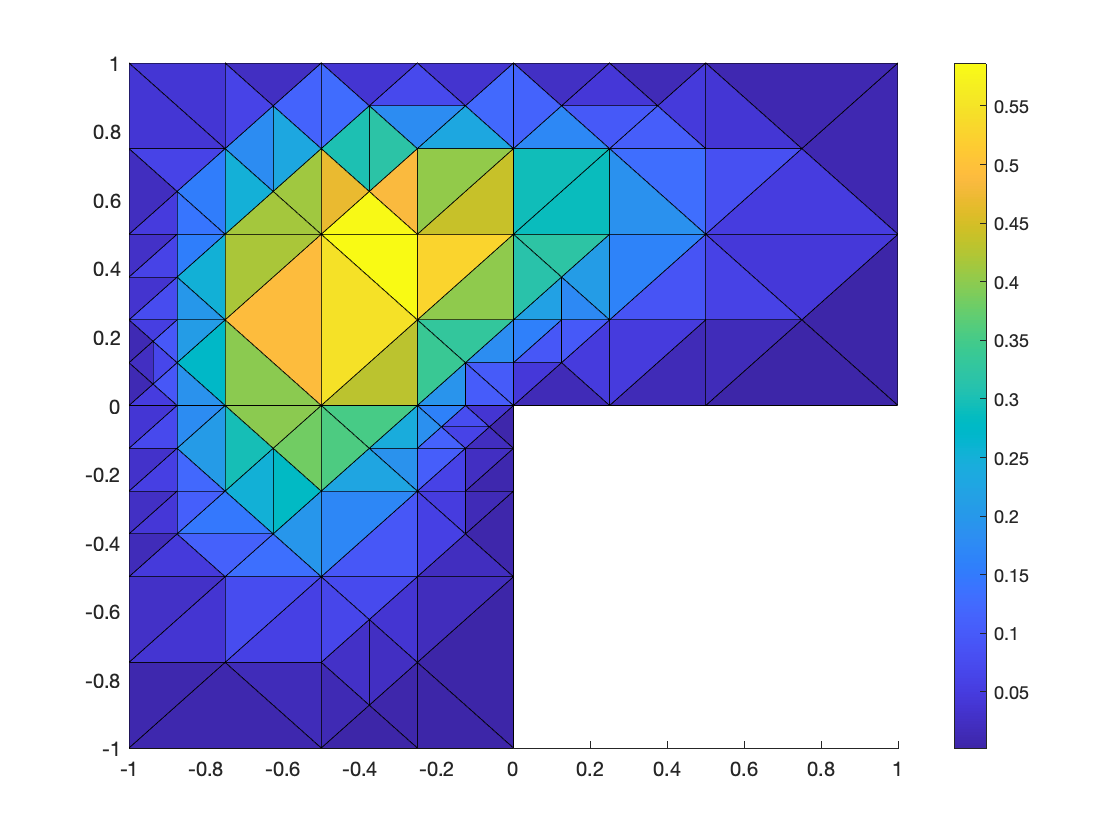}}
		\subfloat[]{\includegraphics[width=0.25\textwidth]{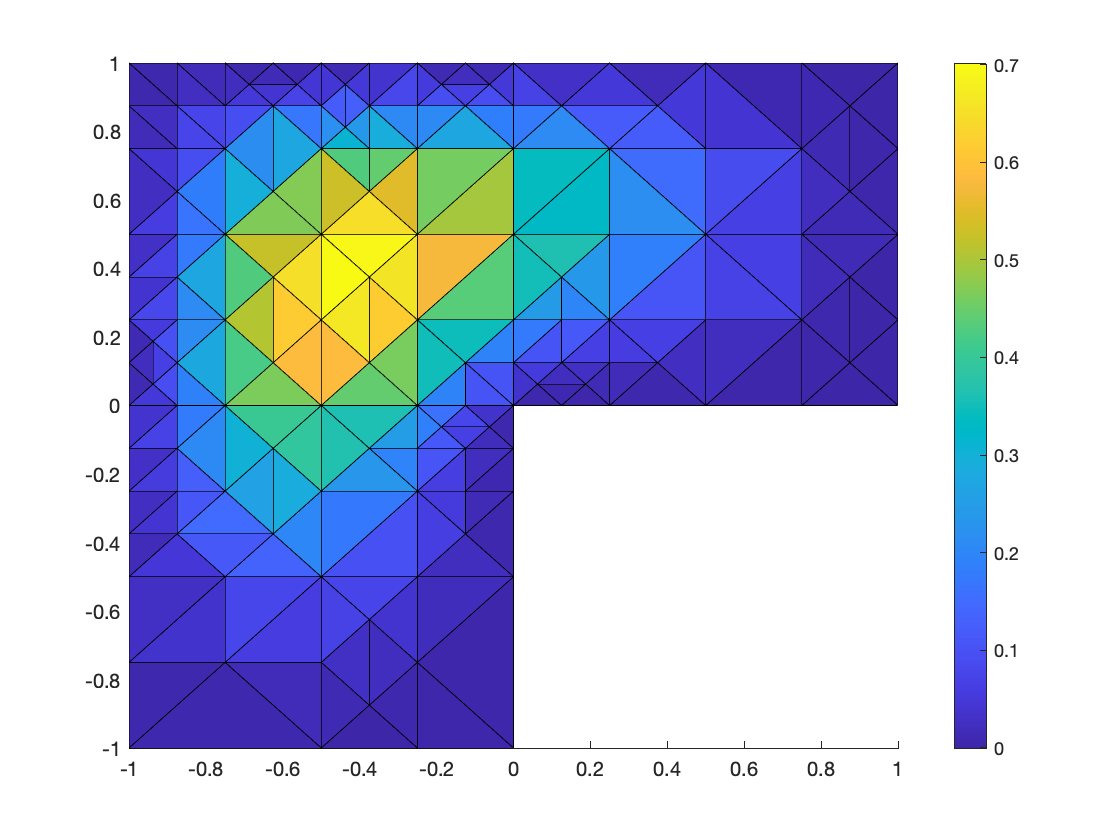}}
		\subfloat[]{\includegraphics[width=0.25\textwidth]{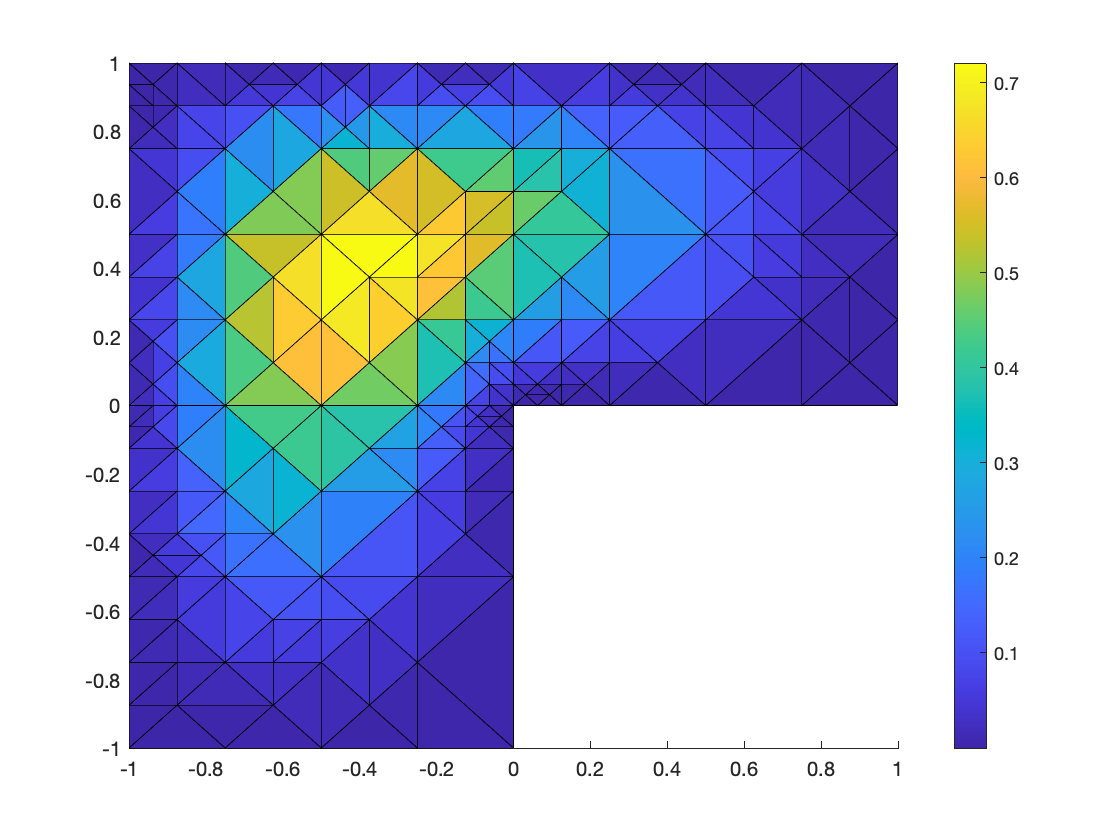}}
		\subfloat[]{\includegraphics[width=0.25\textwidth]{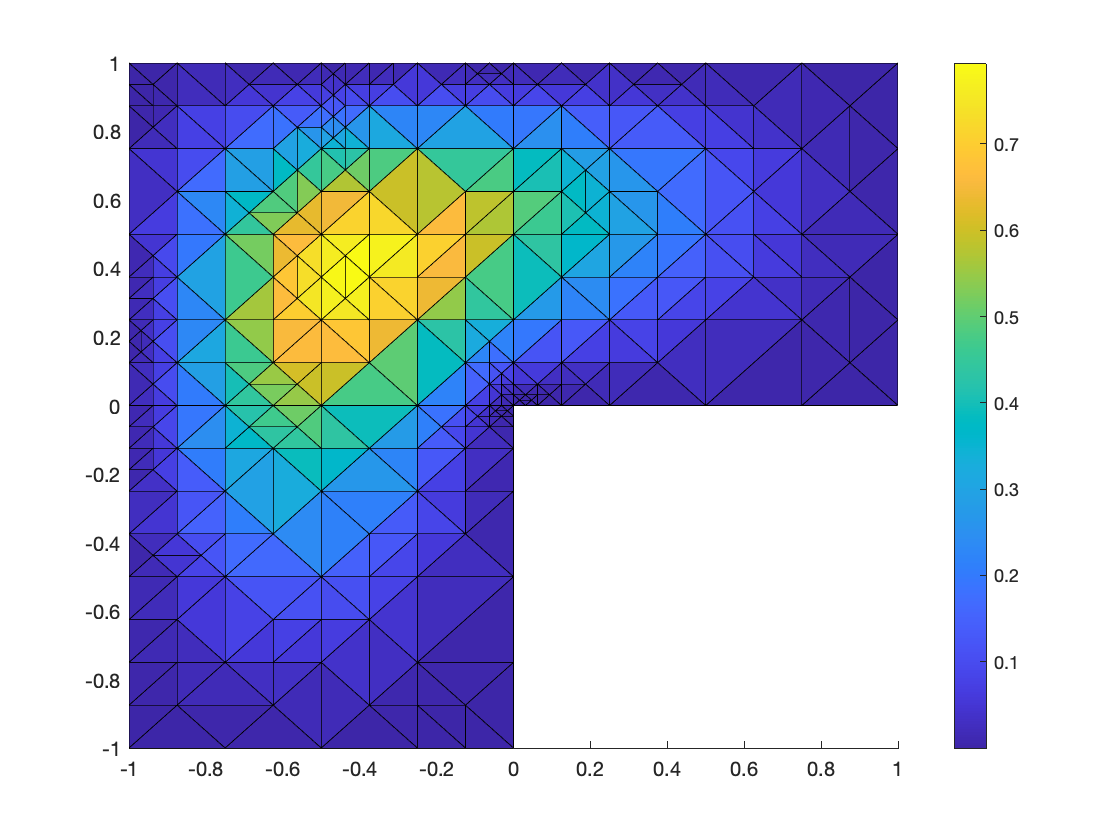}}
	\end{center}
	\caption{ Approximate solution $u_h$ on $\cT_0,\cT_1,\ldots\cT_7$ with parameter $\theta=0.25$ of Algorithm~\ref{alg:adap} for Example~\ref{example_2}}
	\label{fig:uh_Lshape}
\end{figure}

\begin{figure}
	\begin{center}
		\subfloat[]{\includegraphics[width=0.25\textwidth]{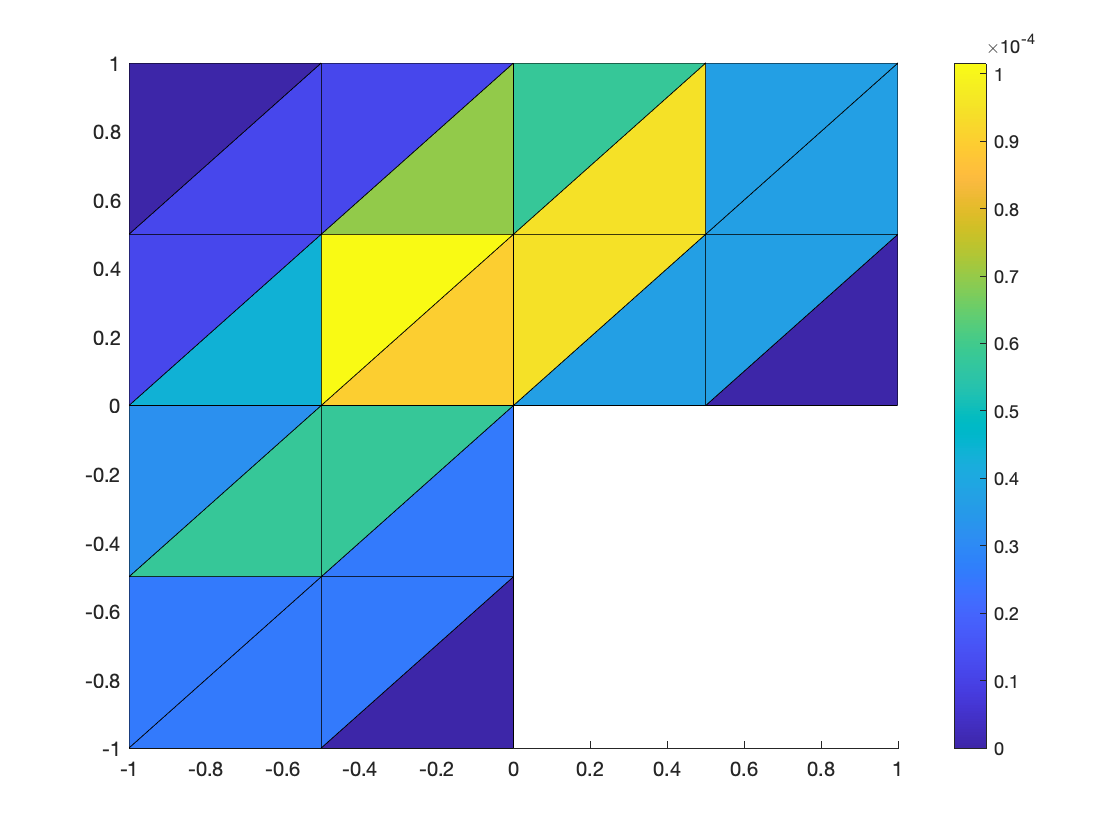}}
		\subfloat[]{\includegraphics[width=0.25\textwidth]{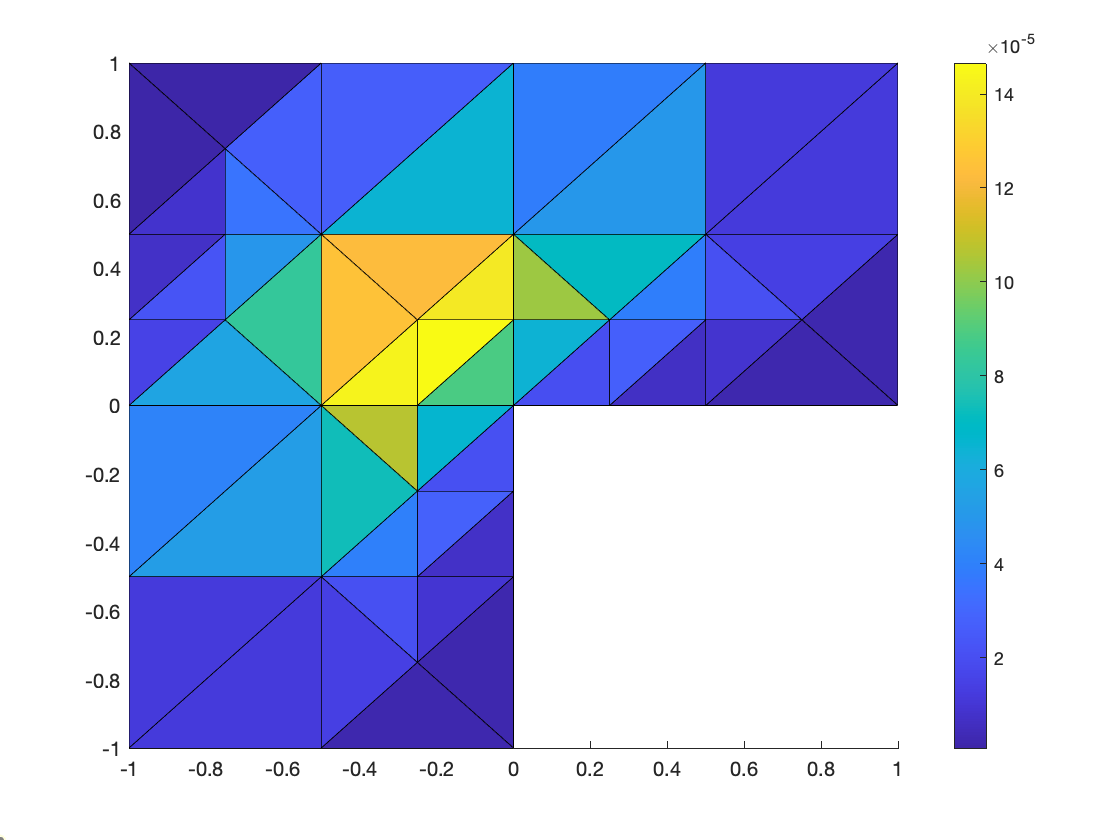}}
		\subfloat[]{\includegraphics[width=0.25\textwidth]{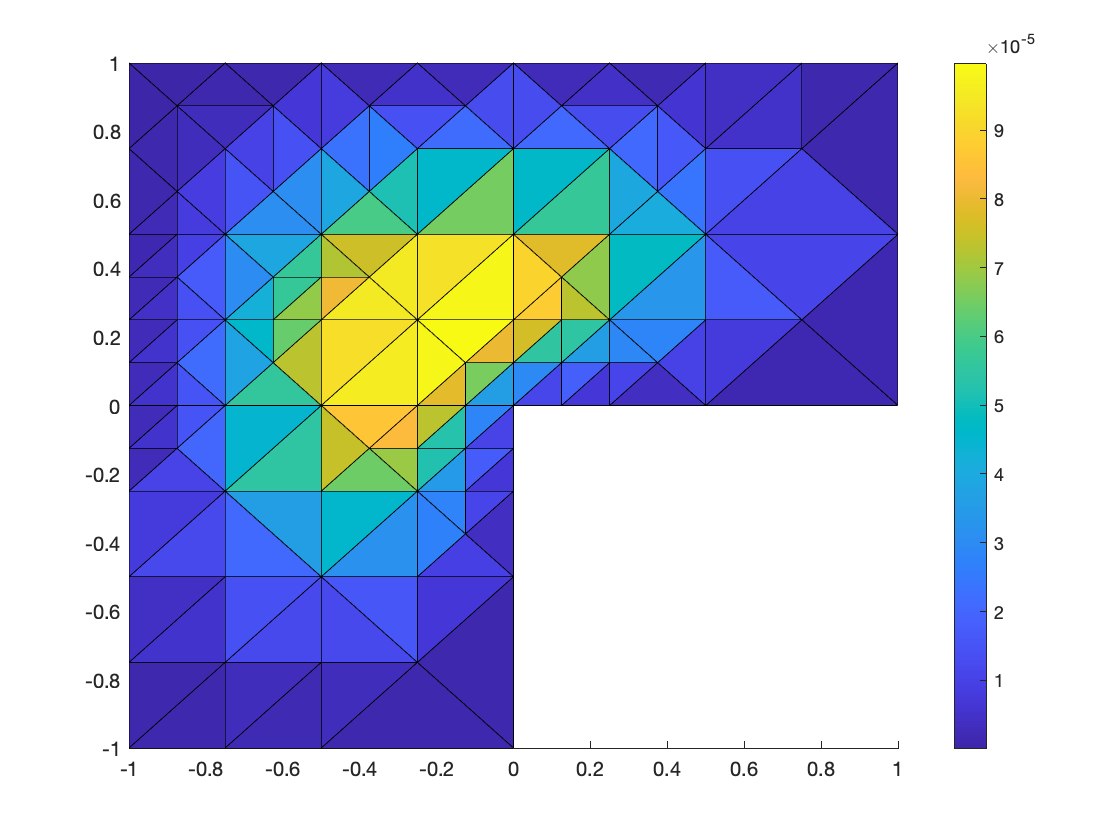}}
		\subfloat[]{\includegraphics[width=0.25\textwidth]{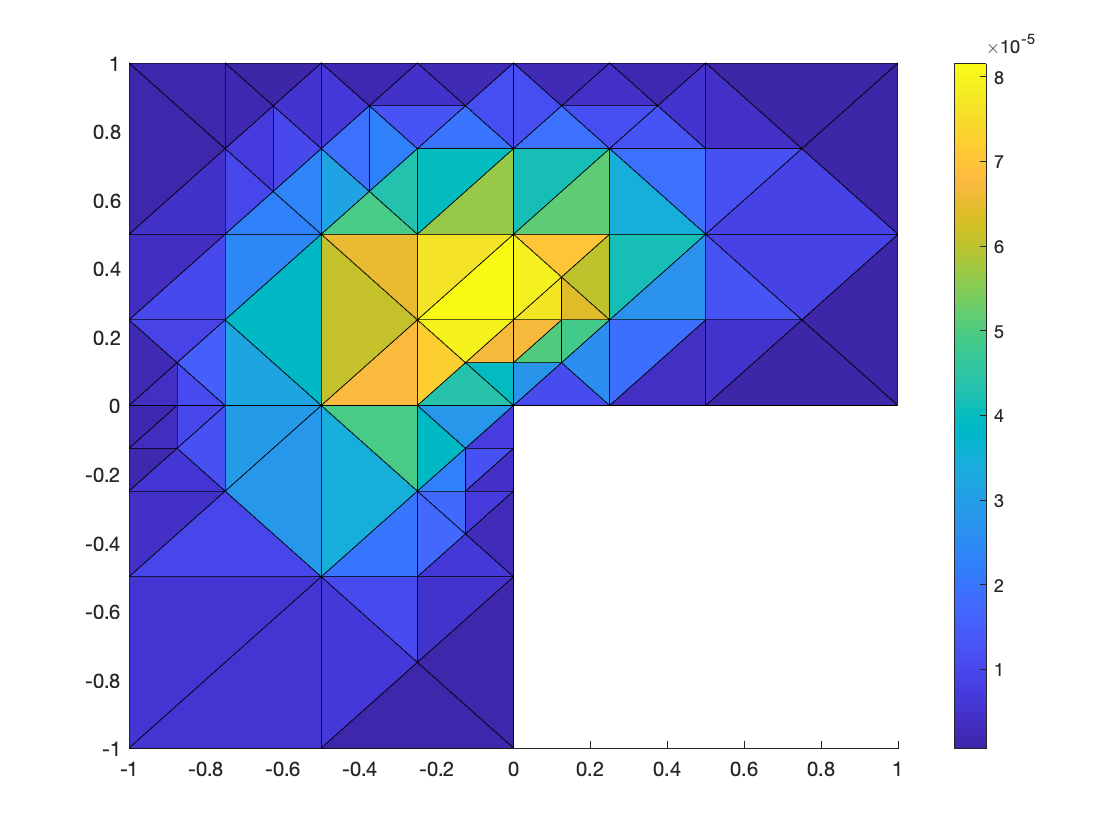}}\\
		\subfloat[]{\includegraphics[width=0.25\textwidth]{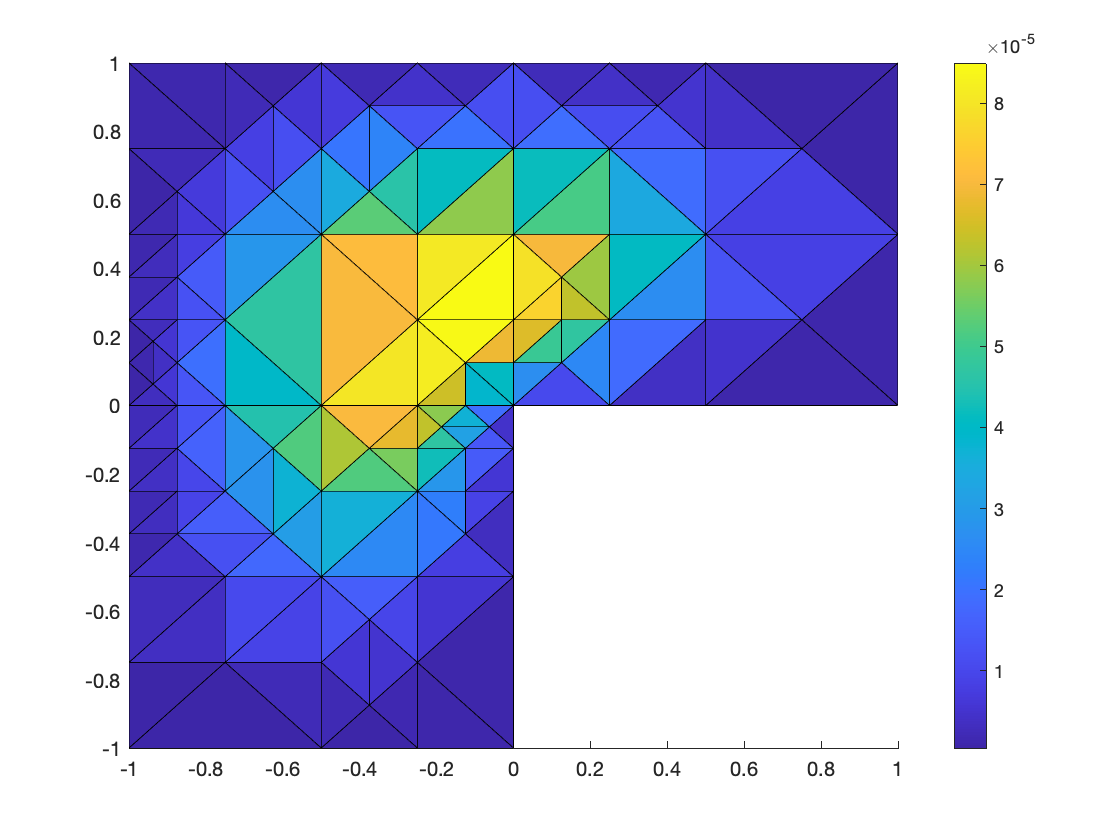}}
		\subfloat[]{\includegraphics[width=0.25\textwidth]{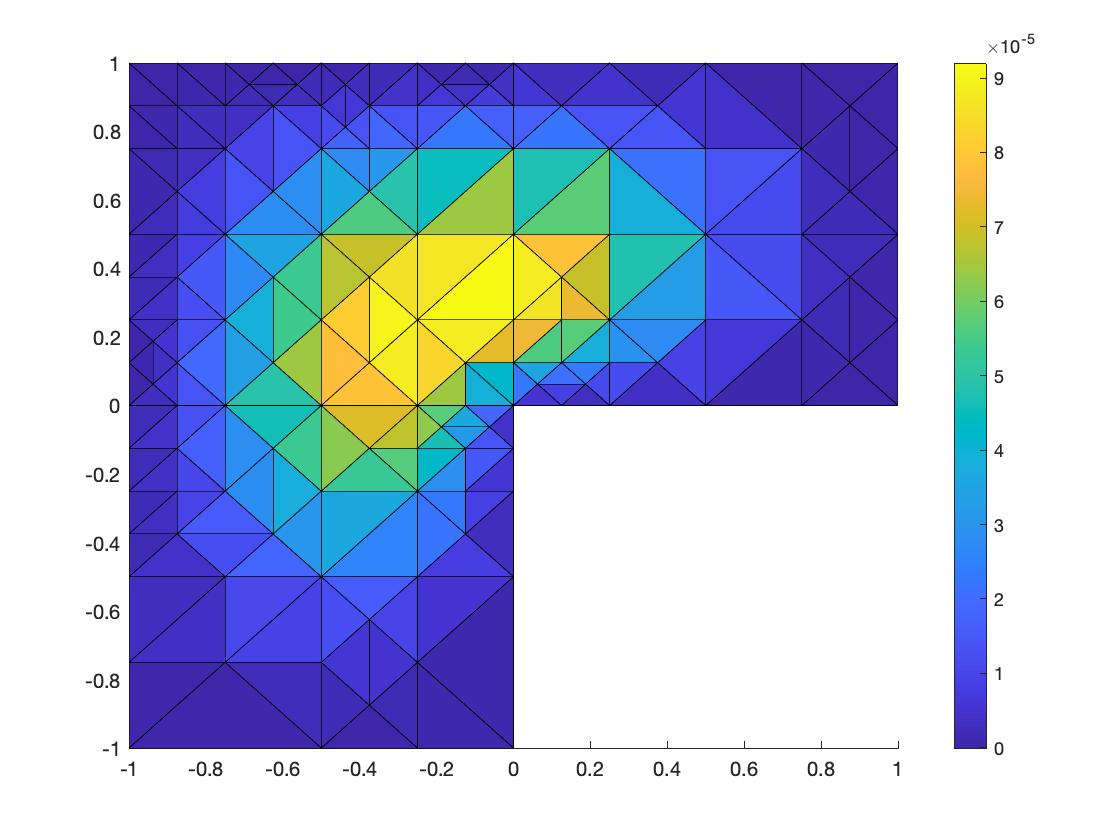}}
		\subfloat[]{\includegraphics[width=0.25\textwidth]{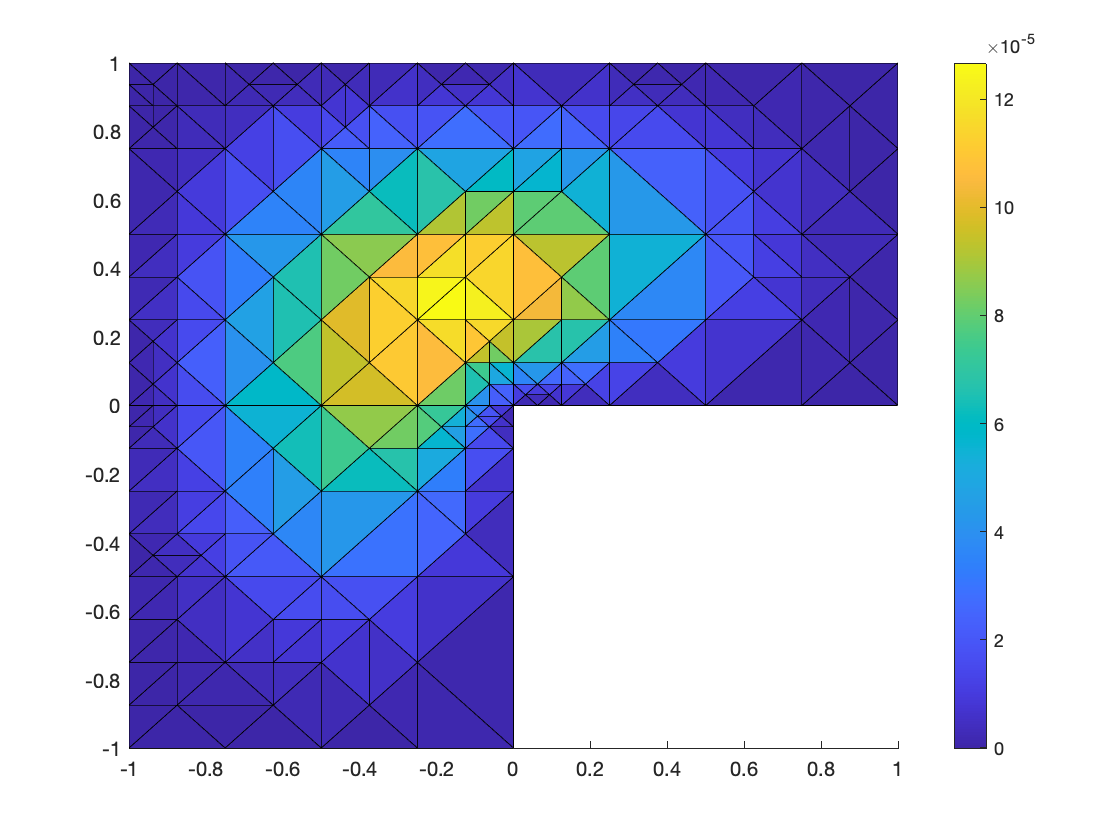}}
		\subfloat[]{\includegraphics[width=0.25\textwidth]{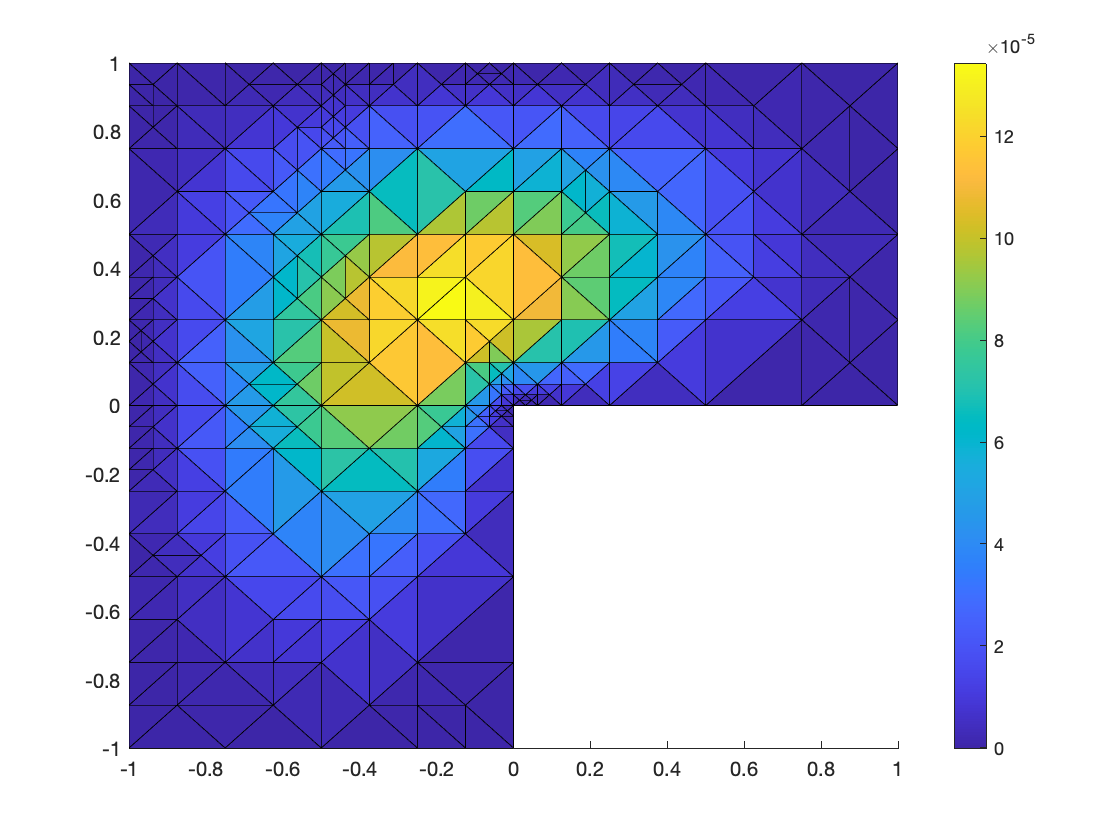}}
	\end{center}
	\caption{ Approximate solution $\tu_h$ on $\cT_0,\cT_1,\ldots\cT_7$ with parameter $\theta=0.25$ of Algorithm~\ref{alg:adap} for Example~\ref{example_2}}
	\label{fig:vh_Lshape}
\end{figure}

\begin{figure}
	\begin{center}
		\subfloat[]{\includegraphics[width=0.4\textwidth]{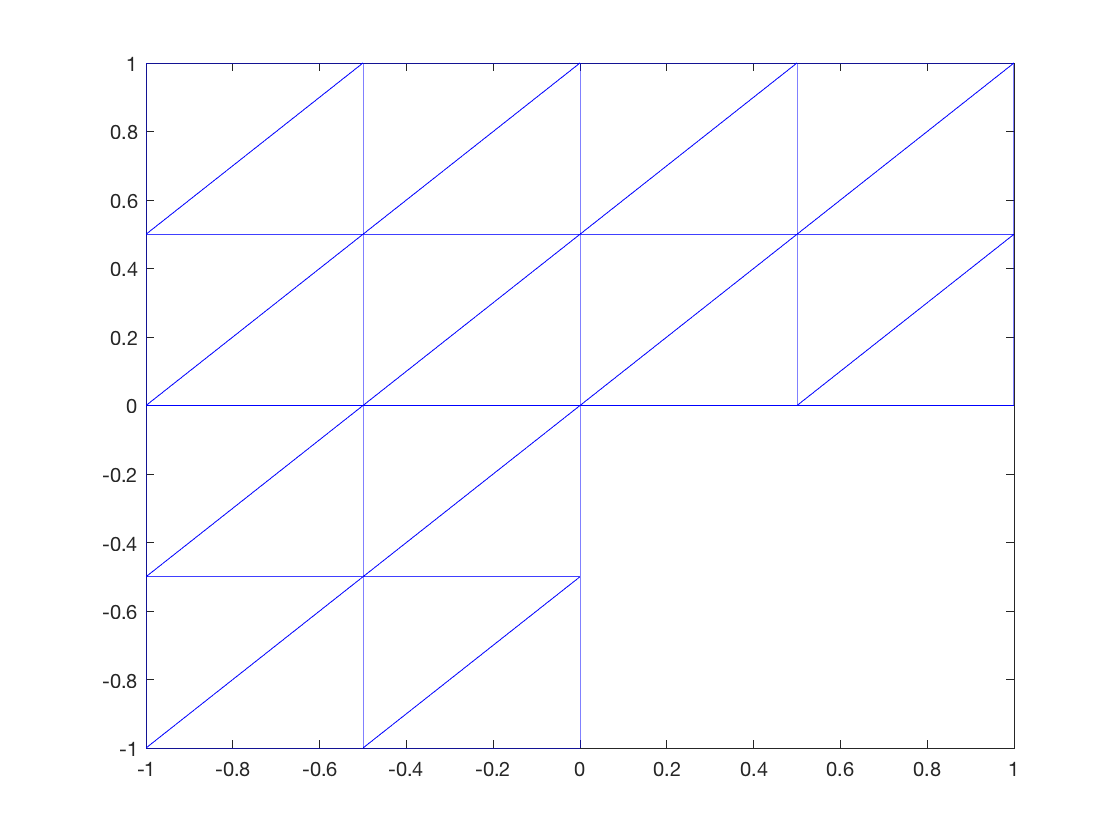}}
		\subfloat[]{\includegraphics[width=0.4\textwidth]{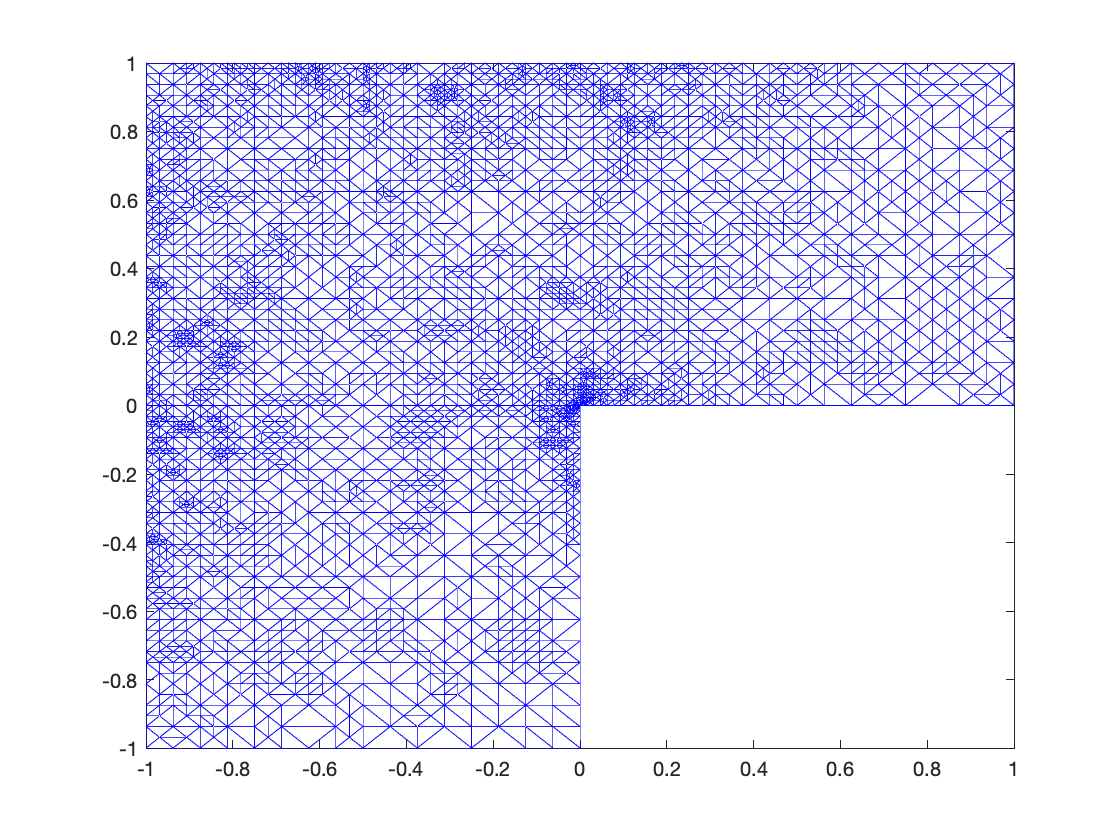}}
	\end{center}
	\caption{(a) Initial triangulation $\cT_0$ and (b) adaptive mesh $\cT_{13}$ with parameter $\theta=0.25$ of Algorithm~\ref{alg:adap} for Example~\ref{example_2}.}
	\label{fig:T_Lshaped}
\end{figure}

\begin{figure}
	\begin{center}
		\includegraphics[width=0.8\textwidth]{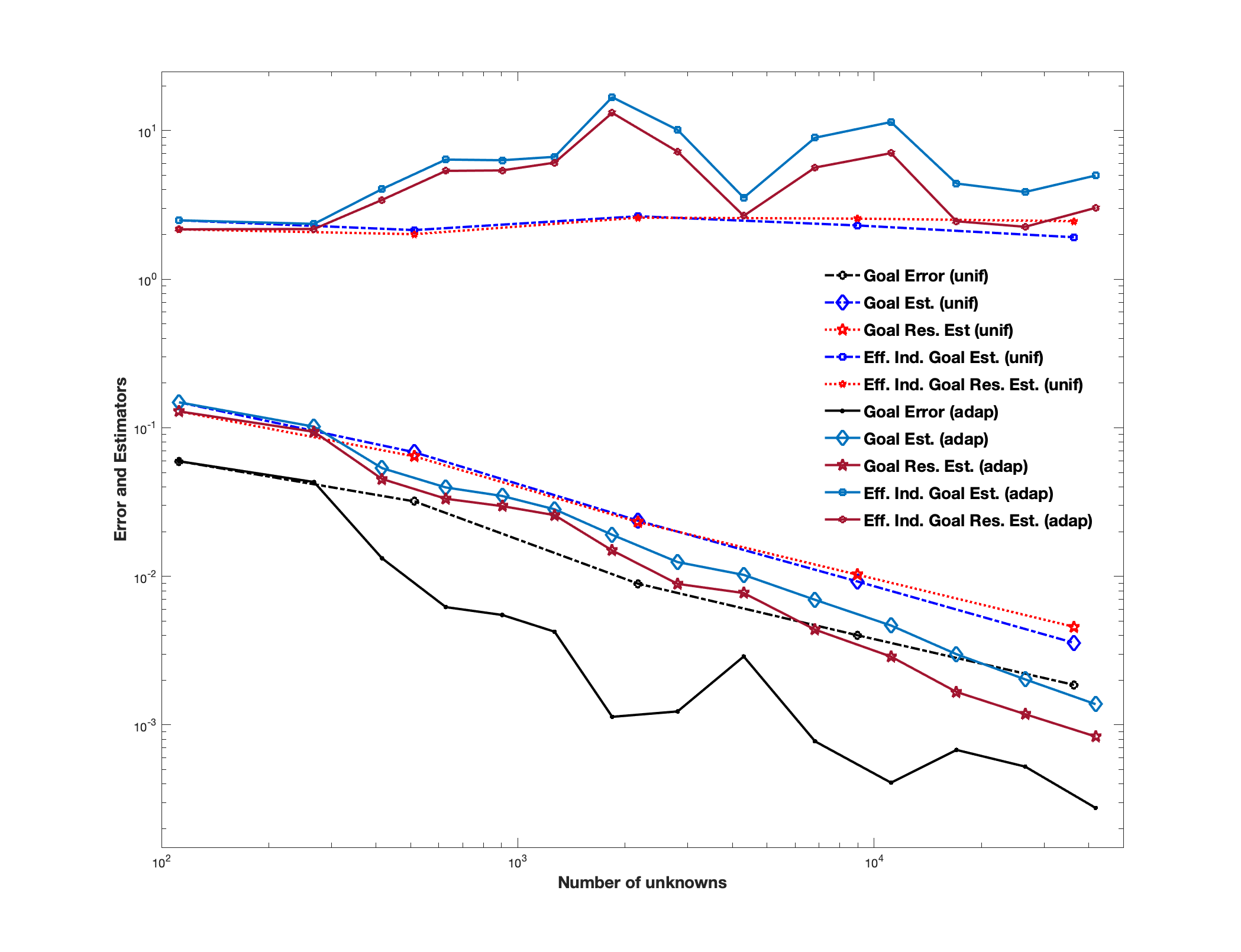}
	\end{center}
	\caption{Convergence histories of goal error and goal estimator with effectivity index for Example~\ref{example_2}.}
	\label{fig:conv_hist_Lshaped}
\end{figure}

For the numerical experiment, we start with an initial mesh $\cT_0$ (see Figure~\ref{fig:T_Lshaped}(a)). We apply the adaptive Algorithm~\ref{alg:adap} with refinement parameter $\theta=0.25$ and maximum refinement level $J=13$ to generate the adaptive meshes $\cT_1,\cT_2,\ldots,\cT_{13}$. We also compare the results with uniform refinement levels $\cT_0,\cT_1,\ldots,\cT_5$. For uniform refinement, each triangle is divided into four similar triangles to obtain next level mesh as described for the previous test. The initial mesh and final adaptive mesh are shown in Figure~\ref{fig:T_Lshaped}. The adaptive meshes and projected solutions for primal and dual problems are illustrated in Figures~\ref{fig:uh_Lshape} \& \ref{fig:vh_Lshape} for first $\cT_0,\cT_1,\ldots,\cT_7$ adaptive meshes. The convergence histories for goal error, goal estimator of \eqref{num_goal_err_est} and \eqref{goal_est_res} with respect to number of unknowns are plotted in Figure~\ref{fig:conv_hist_Lshaped} for uniform and adaptive refinements. It can be observed that for both the refinement procedures goal error is reduced when meshes are refined accordingly. Moreover, the convergence rate for the adaptive refinements is much higher (approximately 2) than the convergence rate for uniform refinement (approximately 1.1). The adaptive algorithm helps to achieve the higher accuracy for the approximation of goal functional with less number of unknowns in the computational process. The effectivity indices of goal estimator $\eta_{h,{\rm goal}}^{\rm abs}$ and goal residual estimator $\eta_{h,{
		\rm res}}$ for uniform refinement appear to be close to $2$ and $2.5$ respectively. Whereas effectivity indices for these estimators for adaptive refinements appear to be close to $5$ and $3$ respectively.

\section{Conclusion}
This article presents an abstract framework for guaranteed goal-oriented a posteriori error control for two variants of discontinuous Galerkin finite element approximation of the model problem~\eqref{eqn:bih}. Also, goal error is represented by an estimator and by a remainder term that combines the dual-weighted residual method and equilibrated moment tensor. The estimators are based on potential reconstruction and equilibrated moment tensor that can be applied to various other finite element approximations for the model problem. The numerical results illustrate that the error in the goal functional \eqref{goal_func} can be controlled efficiently by the Algorithm~\ref{alg:adap}. The methodology described in this article for goal-oriented a posteriori can be applied to nonlinear fourth-order plate problems.

\bigskip
\noindent{\bf Acknowledgements}\\
The author would like to thank  DST C.V. Raman grant R(IA)/CVR-PDF/2020 and National Board for Higher Mathematics (NBHM) research grant 0204/58/2018/R\&D-II/14746 for the financial supports.

\bibliographystyle{siam}
\bibliography{MyBiblio}

\end{document}